%% file: 20260302-rivised.tex
\documentclass[10pt]{amsart}

\usepackage{amsfonts,amsthm,amssymb,amsmath,amscd,ascmac}
\usepackage{latexsym}
\usepackage{lscape}
\usepackage{here}
\usepackage{caption}

\numberwithin{equation}{section}
\newtheorem{lem}{Lemma}[section]
\newtheorem{prop}[lem]{Proposition}
 
\newtheorem{rem}[lem]{Remark} 
\newtheorem{tm}[lem]{Theorem} 

\theoremstyle{definition}
\newtheorem{Def}[lem]{Definition} 
\newtheorem{ex}[lem]{Example} 
 
\newcommand{\NN}{\mathbb{N}} \newcommand{\ZZ}{\mathbb{Z}}
\newcommand{\QQ}{\mathbb{Q}} 
\newcommand{\CC}{\mathbb{C}} \newcommand{\PP}{\mathbb{P}}
\newcommand{\FF}{\mathbb{F}}

\newcommand{\z}{\zeta} \newcommand{\w}{\omega}
\newcommand{\bw}{\overline{\omega}}
\newcommand{\ba}{\overline{\alpha}} \newcommand{\al}{\alpha}
\newcommand{\bb}{\overline{\beta}} \newcommand{\be}{\beta}
\newcommand{\bg}{\overline{\gamma}} \newcommand{\ga}{\gamma}
\newcommand{\bd}{\overline{\delta}} \newcommand{\de}{\delta}

\newcommand{\GL}{\mathrm{GL}} \newcommand{\PGL}{\mathrm{PGL}}
 \newcommand{\PSL}{\mathrm{PSL}}

\newcommand{\Aut}{\mathrm{Aut}} \newcommand{\Out}{\mathrm{Out}}
\newcommand{\Tr}{\mathrm{Tr}} 
\newcommand{\1}{\mathbf{1}} \newcommand{\0}{\mathbf{0}}

\newcommand{\sym}{\mathrm{Sym}} 
\newcommand{\diag}{\mathrm{diag}}

\newcommand{\Q}{\mathrm{Q}} \newcommand{\C}{\mathrm{C}}
\newcommand{\D}{\mathrm{D}} \newcommand{\Sy}{\mathrm{S}}
\newcommand{\A}{\mathrm{A}} \newcommand{\M}{\mathrm{M}}
\begin{document}
\title{Cubic fourfolds with symplectic automorphisms}
\author{Kenji Koike}  
\address{Faculty of Education, University of Yamanashi,
Takeda 4-4-37, Kofu, Yamanashi 400-8510, Japan}
\email{kkoike@yamanashi.ac.jp}
\maketitle
\begin{abstract}
We determine projective equations of smooth complex cubic fourfolds with symplectic 
automorphisms by classifying 6-dimensional projective representations of 
Laza and Zheng's 34 groups. In particular, we determine the number of irreducible 
components for moduli spaces of cubic fourfolds with symplectic actions by these groups. 
We also discuss the fields of definition of cubic fourfolds in six maximal cases.  
\end{abstract}
\input preliminary.tex
\input SmallGroup.tex
\input r16.tex
\input r17.tex
\input r18.tex
\input r19.tex
\input maximal.tex
\appendix 
\renewcommand{\thetable}{\arabic{table}}
\section{An outer automorphism of $\Sy_6$} \label{SS-B}
It is known that automorphisms of $\A_n$ are given by conjugation by elements of 
$\Sy_n$ and $\Aut(\A_n) = \Sy_n$ if $n \ne 6$ and $n \geq 4$. However we have 
$\Aut(\A_6) = \Sy_6.2 = \A_6.2^2$ and it contains three subgroups of index 2, 
that are isomorphic to $\Sy_6, \ \M_{10}$ and $\PGL_2(\FF_9)$ 
(we have $\A_6 \cong \mathrm{PSL}_2(\FF_9)$). Here we give an outer automorphism of 
$\Sy_6$ after \cite{JR82} and related facts.
\\ \indent 
Let us define $\sigma : \Sy_6 \rightarrow \Sy_6$ by
\begin{align*}
\sigma((12)) = (14)(25)(36), \quad \sigma((13)) = (16)(24)(35), 
\quad \sigma((14)) = (13)(26)(45), 
\\
\sigma((15)) = (15)(23)(46), \quad \sigma((16)) = (12)(34)(56).
\end{align*} 
This is a well-defined automorphism of $\Sy_6$ and $\sigma^2$ is a conjugation 
by $(23)(46)$. The automorphism $\sigma$ is an outer automorphism since 
\begin{itemize}
\item the conjugacy class of $(12)$ is swiched with that of $(12)(34)(56)$,
\item the conjugacy class of $(123)$ is swiched with that of $(123)(456)$.
\end{itemize}
We have $\sigma(\A_6) = \A_6$ and $\Aut(\A_6)/\Sy_6 \cong \C_2$ is generated 
by $\sigma$. Moreover we have
\begin{align*}
\sigma((12)(45)) = (13)(46), \quad \sigma((123)) = (123)(465), \quad
\sigma((456)) = (132)(465)
\end{align*} 
and hence $\sigma(\A_{3,3}) = \A_{3,3}$ where 
$\A_{3,3} = (\Sy_3 \times \Sy_3) \cap \A_6$ is generated by
\[
(12)(45), \qquad (123), \qquad (456).
\]
Therefore $\sigma$ gives an outer automorphism of $\A_{3,3}$ 
(see section \ref{SS-A33}). Next let us consider $\Sy_5 \subset \Sy_6$ generated by
$x = (12345)$ and $y = (12)$. Then $\sigma(\Sy_5)$ is generated by
\[
\sigma(x) = (14325), \quad \sigma(y) = (14)(25)(36)
\]
and acts on $\{1,2,3,4,5,6\}$ transitively. Therefore $\sigma(\Sy_5)$ is not
conjugate to a standard $\Sy_5 \subset \Sy_6$. Taking intersections with $\A_6$, 
we obtain two types of $\A_5$ in $\A_6$. For this reason, we have two families of 
cubic fourfolds with symplectic actions of $\A_5$ (see section \ref{SS-A5}). 
\section{Table}
For non-maximal cases of Laza and Zheng's 34 groups, we have the following table, 
where $r(G) = \mathrm{rank} \ S_G$ (see \cite{LZ22}), GAP-ID is ``IdGroup($G$)'' in GAP, 
$M$ is the Schur multiplier of $G$, and ID* is GAP-ID of a Schur covering group $G^*$ 
given by ``SchurCover($G$)'' in GAP. Equations of corresponding cubic 4-folds are given in Eq, 
and irr is the number of irreducible components of the moduli space $\mathcal{M}(G)$ 
of cubic fourfolds $X$ such that $G \subset \Aut^s(X)$. 
\begin{table}[hbtp] 
\caption{Schur covers and multipliers} \label{TBL-Schur}
\begin{tabular}{c|c|c|c|c|c|c|c} 
\hline
 $r(G)$ & $G$ & GAP-ID & $M$ & ID* & Eq & irr & \\ 
\hline
8 & 2 & 2,1 & 1 & 2,1 & (\ref{eqC2}) & 1 \\
12 & $2^2$ & 4,2 & 2 & 8,3 & (\ref{eqC2xC2}) & 1 \\ 
12 & 3 & 3,1 & 1 & 3,1 & (\ref{eqC3-1}), (\ref{eqC3-2}) & 2 \\
14 & 4 & 4,1 & 1 & 4,1 & (\ref{eqC4}) & 1 \\ 
14 & $\Sy_3$ & 6,1 & 1 & 6,1 & (\ref{eqS3}) & 1\\
15 & $\D_8$ & 8,3 & 2 & 16,7 & (\ref{eqD8}) & 1\\
16 & $\A_{3,3}$ & 18,4 & 3 & 54,8 & (\ref{eqA33-1}), (\ref{eqA33-2}) & 2 \\
16 & $\D_{12}$ & 12,4 & 2 & 24,8 & (\ref{eqD12-1}), (\ref{eqD12-2}) & 2 \\ 
16 & $\A_4$ & 12,3 & 2 & 24,3 & (\ref{eqA4-1}), (\ref{eqA4-2}) & 2 \\ 
16 & $\D_{10}$ & 10, 1 & 1 & 10, 1 & (\ref{eqD10}) & 1 \\
17 & $\Sy_4$ & 24,12 & 2 & 48,28 &  (\ref{eqS4-1}), (\ref{eqS4-2})  & 1 \\ 
17 & $\Q_8$ & 8,4 & 1 & 8,4 & (\ref{eqQ8}) & 1 \\
18 & $3^{1+4}_+:2$ & 486, 249 & $3^5$ & -& (\ref{eq3^{1+4}:2}) & 1 \\ 
18 & $\A_{4,3}$ & 72,43 & 6 & 432,258 & (\ref{eqA43-1}), (\ref{eqA43-2}) & 2 \\
18 & $\A_5$ & 60, 5 & 2 & 120,5 & (\ref{eqA5-1}), (\ref{eqA5-2}) & 2 \\ 
18 & $3^2.4$ & 36,9 & 3 & 108,15 & (\ref{eq3^2.4-1}), (\ref{eq3^2.4-2}) & 2 \\ 
18 & $\Sy_{3,3}$ & 36,10 & 2 & 72,22 & (\ref{eqS33}) & 1 \\ 
18 & $\mathrm{F}_{21}$ & 21,1 & 1 & 21,1 & (\ref{eqF21}) & 1 \\
18 & $\mathrm{Hol}(5)$ & 20, 3 & 1 & 20,3 & (\ref{eqHol(5)})  & 1 \\ 
18 & $\mathrm{QD}_{16}$ & 16,8 & 1 & 16,8 & (\ref{eqQD16}) & 1 \\
19 & $3^{1+4}_+:2.2$ & 972, 776 & $3^4$ & - & (\ref{eq3^{1+4}:2.2}) & 1 \\ 
19 & $\A_6$ & 360,118 & 6 & - & (\ref{eqA6-1}),  (\ref{eqA6-2}) & 2 \\
19 & $\mathrm{L}_2(7)$ & 168,42 & 2 & 336,114 & (\ref{eqL2(7)}) & 1 \\ 
19 & $\Sy_5$ & 120,34 & 2 & 240, 89 & (\ref{eqS5-1}), (\ref{eqS5-2}) & 2 \\ 
19 & $\M_9$ & 72,41 & 3 & 216,88 & (\ref{eqM9}) & 1 \\ 
19 & $\mathrm{N}_{72}$ & 72,40 & 2 & 144,117 & (\ref{eqN72}) & 1 \\ 
19 & $\mathrm{T}_{48}$ & 48,29 & 1 & 48,29 & (\ref{eqT48}) & 1 \\
\hline
\end{tabular}
\end{table}

\end{document}

%% file: preliminary.tex
\section{Introduction}
In \cite{LZ22}, Laza and Zheng completely classified the groups $\Aut^s(X)$ 
of symplectic automorphisms of complex cubic fourfolds $X$.  
According to them, we have 34 possible groups as $\Aut^s(X)$ 
with six maximal cases
\begin{align*}
3^4:\A_6, \quad \A_7, \quad 3^{1+4}_+:2.2^2, \quad \M_{10}, \quad \mathrm{L}_2(11), \quad \A_{3,5}.
\end{align*}
 For each group $G$ among the 34 groups, we have the covariant lattice $S_G$ defined as the 
orthogonal complement of the invariant lattice $H^4(X,\ZZ)^G$. 
It is uniquely determined (up to isomorphism) by $G$, and the moduli space of cubic fourfolds $X$ 
such that $G \subset \Aut^s(X)$ has dimension $20 - \mathrm{rank} \, S_G$.
However, this moduli space is not necessarily irreducible, corresponding to the fact that there may exist 
different primitive embeddings $S_G \hookrightarrow H^4(X,\ZZ)$. 
For example, there exists the unique cubic fourfold with $\Aut^s(X) = G$ for the above maximal 
groups $G$ except $\A_7$ and $\M_{10}$, but two different cubic fourfolds if 
$\Aut^s(X) = \A_7$ or $\M_{10}$.
Laza and Zheng classified possible finite groups $G = \Aut^s(X)$ via embeddings 
$(G, S_G) \hookrightarrow (\mathrm{Co}_0, \mathbb{L})$, where $\mathbb{L}$ 
is the Leech lattice and $\mathrm{Co}_0 = \mathrm{O}(\mathrm{L})$ is the Conway group, 
and the result is based on the classification of the fixed-point lattices for the Leech lattice by 
H\"ohn and Mason (\cite{HM16}). 
\\ \indent
More recently, Yang, Yu and Zhu studied full automorphism groups for cubic fivefolds and fourfolds 
in \cite{YYZ23}. According to them, there are 15 maximal groups, and a finite group $G$ can acts 
faithfully on a smooth cubic fourfold if and only if $G$ is a subgroup of these 15 groups. 
However, explicit equations of cubic fourfolds with automorphisms are still unknown in many cases.  
\\ \indent
In this paper, we study symplectic actions of 34 groups of Laza and Zheng, and detertmine invariant 
smooth cubic fourfolds. Since any automorphism of a smooth cubic fourfold is given by a projective 
transformation, our problem is reduced to classification of faithfull projective representations 
of 34 groups which satisfy symplectic conditions. 
We classify such projective representations by considering character tables of Schur covers 
for 34 groups. A Schur cover $G^*$ of a finite group $G$ is defined as a central extension $G^* = M.G$ 
such that $M \cong \mathrm{H}^2(G, \CC^{\times})$ and $M \subset [G^*, G^*]$. 
It is known that complex projective representations of a finite group $G$ lift to linear representations 
of a Schur cover $G^*$. Therefore, for a fixed $G$, we can obtain all families of cubic fourfolds 
such that $G \subset \Aut^s(X)$ in the following manner. 
\\ \\
(1) \ Determine linear representations $\rho : G^* \rightarrow \GL_6(\CC)$ which induces a 
faithful projective action of $G$ and satisfies symplectic conditions. 
\\
(2) \ Check that the vector space of $\rho(G)$-invariant cubic forms $V_3(\rho) \subset \CC[x_1, \cdots, x_6]$  
contains $F(x) \in V_3(\rho)$ which gives  a smooth cubic fourfold $X_F$. 
\\
(3) \ Check that $\Aut^s(X_F) = G$ for a general member $F \in V_3(\rho)$.  
\\ \\ 
For (1),  we can determine possible representations from character tables relying on 
absolute values of characters (see section \ref{symplectic}). We can check (3) by counting 
dimensions of moduli spaces (see section \ref{ModuliSpace}). 
For each representation $\rho$ obtained in this manner, $V_3(\rho)$ defines an irreducible 
family of cubic fourfolds such that $G \subset \Aut^s(X)$. However, we may have 
$V_3(\rho_1) = V_3(\rho_2)$ even if $\rho_1$ is not projectively equivalent to $\rho_2$. 
This is caused by the action of $\Out(G)$. For an outer automorphism $\sigma$ of $G$, 
we may have $\rho \ne \rho \circ \sigma$ although we have $V_3(\rho) = V_3(\rho \circ \sigma)$. Excluding duplicated families, we have
\begin{tm} 
The moduli space of cubic fourfolds with $G \subset \Aut^s(X)$ has two irreducible components 
if $G$ is one of the following groups:
\begin{align*}
3, \ \A_{3,3}, \ \D_{12}, \ \A_4, \ \A_{4,3}, \ \A_5, \ 3^2.4, \ \A_6, \ \Sy_5, 
\ \A_7, \ \M_{10}.
\end{align*}
For other cases, the moduli spaces are irreducible {\rm (}see Table \ref{TBL-Schur} in 
Appendix B \rm{)}.
\end{tm}
For maximal cases, defining equations of corresponding cubic fourfolds are known, 
and they are defined over $\QQ$ in most cases. In the last section, we show
\begin{tm} 
{\rm (1)} \ If the symplectic automorphisms $\Aut^s(X)$ of a cubic fourfold $X$ is maximal, then $X$ is 
defined over $\QQ$ except in the case of $\Aut^s(X) = \M_{10}$. 
\\
{\rm (2)} \ Let $X_+$ and $X_-$ be two different cubic fourfolds 
with $\Aut^s(X_{\pm}) = \M_{10}$. They are defined over 
$K = \QQ(\sqrt{6})$, and conjugate to each other. Namely, $\mathrm{Gal}(K/\QQ)$ 
swiches $X_+$ and $X_-$.  
\end{tm}
Our results are heavily dependent on computer system GAP. 
Also we used databases \cite{AtlasWeb} and \cite{GroupNames} for information of 
finite groups such as generators and automorphisms. 
\\ \\
{\bf Acknowledgements.} 
The author is grateful to referees for many suggestions to improve the manuscript.
It was raised by referees that the number of irreducible components is  
(perhaps more easily) computed by computing the number of primitive 
embeddings of the lattice $S_G = (H^4(X, \ZZ)_{prim})_G$ into $H^4(X, \ZZ)_{prim}$. 
The number of embeddings is included in \cite{MM24}, Table 4. 
The groups of Laza and Zheng are included in this list, and each time they appear
with a $\times$ in the cubic column corresponds to the number of
primitive embeddings of $S_G$.
\\ \\
{\bf Notations.} 
We say that $G$ is an {\it extension of $B$ by $A$} and denote $G \cong A.B$
if there is an exact sequence
\begin{align*}
1 \longrightarrow A \overset{\iota}{\longrightarrow} G 
\overset{\pi}\longrightarrow B \longrightarrow 1
\end{align*}
(the structure of $A.B$ is not uniquely determined).  In particular, it is called 
a {\it central extension} if $\iota(A)$ is included in the center of $G$.  We denote 
$G \cong A:B$ if the exact sequence is split, that is, if $\pi$ has a section. 
\\ \indent
We use the following notations:
\begin{itemize}
\item $\C_m, \ \D_m, \ \Q_m$ denote {\it cyclic, dihedral, quaternion groups} of order $m$.
\item $\C_m$ is often denoted by $m$ (e.g. $4$ denotes $\C_4$, and $2^2$ denotes 
$\C_2 \times \C_2$). 
\item $p^{1+2n}_{\pm}$ denote {\it extraspecial groups} of order $p^{1+2n}$ 
(see Remark \ref{extraspecial}).
\item $\Sy_{m,n} = \Sy_m \times \Sy_n, \ \A_{m,n} = \Sy_{m,n} \cap \A_{m+n}$, where 
$\A_n, \ \Sy_n$ denote {\it alternating, symmetric groups} of degree $n$.
\item $\diag(a_1,\dots, a_n)$ denotes a diagonal matrix with diagonal entries $a_1,\dots, a_n$.
\item $A \oplus B$ denotes a block matrix 
$\begin{bmatrix} A & \bf{0} \\ \bf{0} & B \end{bmatrix}$ for square matrices $A$ and $B$
\item $\zeta_n = \exp(2\pi i/n), \quad \omega = \zeta_3$, \quad 
$\widehat{G} = \mathrm{Hom}(G, \CC^{\times})$, \quad$s_k(x) = x_1^k + \cdots + x_n^k$.
\end{itemize}
\section{Projective actions of finite groups and invariant hypersurfaces}
\subsection{Schur covering groups}
Let $X \subset \PP^{n+1}(\CC)$ be a hypersurface of degree $d$. If $X$ is smooth and $n \geq 2, \ d \geq 3$, 
then the group $\Aut(X)$ of automorphisms is finite and realized by projective transformations of 
$\PP^{n+1}$ except the case $(n,d)=(2,4)$ (\cite{MM64}).  
First of all, we recall basic facts on projective representations of finite groups. For details, see \cite{Ka85}. 
\\ \indent
A projective representation of a group $G$ is a homomorphism $\rho : G \rightarrow \PGL_n(\CC)$.  
If there is a commutative diagram 
\[
  \begin{CD}
     \tilde{G} @> \tilde{\rho} >> \GL_{n+1}(\CC)  \\
     @VVV @VV \pi V \\
     G @> \rho >> \PGL_{n+1}(\CC)
  \end{CD}
\]
of group homomorphisms, we say that $\rho$ is lifted to a linear representation $\tilde{\rho}$. 
Let us see actions of the alternating group $\A_7$ as a typical example. 
Note that any non-trivial action of a non-commutative simple group on a smooth cubic fourfold 
is symplectic (see section \ref{symplectic}). 
\begin{ex} \label{exA7}
According to \cite{LZ22}, there exist two non-isomorphic smooth cubic fourfolds with symplectic 
action of $\A_7$. One is the Clebsch-Segre cubic
\begin{align} \label{Clebsch-Segre}
\begin{cases}
x_1 + x_2 + \cdots + x_7 = 0 \\
x_1^3 + x_2^3 + \cdots + x_7^3 = 0
\end{cases} 
\text{in} \ \PP^{6},
\end{align}
and the symmetric group $\Sy_7$ acts on it as permutations of coordinates 
(However, the action of $(12) \in \Sy_7$ is not symplectic).  
Another $\A_7$-invariant cubic fourfold was identified by Yang, Yu and Zhu in \cite{YYZ23} 
considering the ``triple cover'' of $\A_7$. Their model is defined over $\QQ(\sqrt{15})$, and 
we give an equation over $\QQ$ here. 
Let $G$ be a subgroup of $\GL_6(\CC)$ generated by two matrices
\begin{align} \label{ab in 3.A7}
a = \begin{bmatrix}
 1 & 0 & 0 & 0 & 0 & 0 \\
 0 & 1 & 0 & 0 & 0 & 0 \\
 0 & 0 & \w & 0 & 0 & 0 \\
 0 & 0 & 0 & \w & 0 & 0 \\
 0 & 0 & 0 & 0 & \w^2 & 0 \\
 0 & 0 & 0 & 0 & 0 & \w^2 \\
\end{bmatrix},
\quad
b = \frac{1}{3} \begin{bmatrix}
 -1 & 0 & -\w & 0 & 3 \w^2 & 2 \w^2 \\
 2 & 0 & 2 \w & 0 & 0 & -\w^2 \\
 -\w^2 & 0 & 2 & 3 & 0 & -\w \\
 2 \w^2 & 0 & -1 & 0 & 0 & 2 \w \\
 2 \w & 3 \w & -\w^2 & 0 & 0 & -1 \\
 -\w & 0 & 2 \w^2 & 0 & 0 & 2 \\
\end{bmatrix}. 
\end{align}
We can check that $G \cong 3.\A_7$ by GAP. More precisely, the center $\C_3$ is 
generated by $(abab^{-1})^2 = \diag (\w, \dots, \w)$, and we have 
a commutative diagram:  
\[
   \begin{CD}
     \C_3 @>>> G @>>> G/\C_3 \cong \A_7 \\
  @VVV    @VVV @VVV \\
    \CC^{\times}   @>>>  \GL_6(\CC) @>>> \PGL_6(\CC)
  \end{CD}.
\]
The group $G \cong 3.\A_7$ is a non-split extension, and this projective action of 
$\A_7$ cannot be lifted to a linear representation of $\A_7$ itself.  
We have the unique invariant cubic 4-fold $X$:
\begin{multline} \label{3.A7-invariant}
f(x) = 2(x_1^3 + x_2^3 + x_3^3 + x_4^3 + x_5^3 + x_6^3) \\
+ 3 (x_1^2 x_2 + x_1 x_2^2 + x_3^2 x_4 + x_3 x_4^2 + x_5^2 x_6 + x_5 x_6^2) \\
+2 (x_2 x_3 x_5 + x_1 x_4 x_5 + x_1 x_3 x_6 + x_2 x_4 x_6) \\
+4 (x_1 x_3 x_5 + x_2 x_4 x_5 + x_2 x_3 x_6 + x_1 x_4 x_6) = 0.
\end{multline}
Considering $X$ by reduction modulo 2, it is easily checked that $X$ is smooth 
(In fact, $X$ is isomorphic to the Fermat cubic over $\mathbb{F}_4$).  
\hfill $\square$
\end{ex}
In general, projective representations (over $\CC$) of a finite group $G$ can be lifted to linear 
representations of a central extension of $G$, called a Schur covering group. More precisely,  a group 
$G^*$ is called a {\it Schur covering group} (or {\it representation group}) of $G$ if $G^*$ is 
a central extension
\begin{align} \label{SC sequence}
1 \longrightarrow M \overset{\iota}{\longrightarrow} G^* \longrightarrow G \longrightarrow 1
\end{align}
such that $M \cong \mathrm{H}^2(G, \CC^{\times})$ and $\iota(M) \subset [G^*, G^*]$. 
In general, $G^*$ is not uniquely determined. The Klein's four group $\C_2 \times \C_2$ has two 
Schur covering groups $\D_8$ and $\Q_8$. However, it is known that all projective representations of $G$ 
can be lifted to linear representations of one Schur covering group $G^*$. 
The abelian group $M = M(G)$ in equation (\ref{SC sequence}) is called the {\it Schur multiplier} of $G$.  
For example, we have
\begin{align*}
M(\A_n) = \begin{cases} 1 \quad (n \leq 3) \\ \C_6 \quad (n = 6,7) \\
\C_2 \quad (\text{others}) \end{cases}, \quad
M(\Sy_n) = \begin{cases} 1 \quad (n \leq 3) \\ \C_2 \quad (n \geq 4) \end{cases}
\end{align*}
(see \cite{Ka85}). If $M(G) =1$, we have $G^* =G$ and any projective representation of $G$ 
is lifted to a linear representation of $G$ itself. 
For non-maximal cases of Laza and Zheng's 34 groups, Schur covers are given in Appendix B.  
\begin{rem} \label{RM-S-number}
It is known that $G$ has only one Schur covering group if the orders of $G/[G,G]$ and 
$M(G)$ are coprime {\rm (\cite{Ka85}, Theorem 4.4)}. 
\end{rem}
\subsection{Symplectic automorphisms of cubic 4-folds} \label{symplectic}
Let $X$ be a smooth cubic fourfold. Then $\Aut(X)$ acts faithfully on $\mathrm{H}^4(X, \ZZ)$. 
If $g \in \Aut(X)$ acts trivially on $\mathrm{H}^{3,1}(X) \cong \CC$, the action of $g$ is called  
symplectic. Therefore we have an exact sequence 
\[
1 \longrightarrow \Aut^s(X) \longrightarrow \Aut(X) \longrightarrow 
\Aut(\mathrm{H}^{3,1}(X)) = \CC^{\times}. 
\]
for the group $\Aut^s(X)$ of symplectic automorphisms. 
In particular, a non-trivial action of a non-commutative simple group is necessarily symplectic. 
\\ \indent
In \cite{Fu16}, smooth cubic fourfolds with symplectic automorphisms of primary order were classified 
by Fu. The result was generalized to symplectic automorphisms of any order in \cite{LZ22}. 
According to them, the order of a symplectic automorphism is one of 
\[
1, 2, 3, 4, 5, 6, 7, 8, 9, 11, 12, 15.
\]
If a cubic fourfold $X$ has a symplectic automorphism of order $n = 9, 11, 12, 15$, then we have
\[
\Aut^s(X) = 3^4:\A_6, \quad \mathrm{L}_2(11), \quad 3^{1+4}:2.2^2, \quad \A_{3,5}, 
\]
respectively, and the cubic fourfold $X$ is unique for these cases.  
Symplectic automorphisms $g$ of order $n \leq 8$ are classified into the following 
diagonalized normal forms:
\begin{table}[H]
  \caption{Normal forms of symplectic automorphisms} \label{TBL-Fu}
\begin{tabular}{c||c|c}
 $n$ & $(k_1, \cdots, k_6)$ & $|\Tr \, g|$ \\
\hline
2 & (0,0,0,0,1,1) & 2 \\
\hline
3 & (0,0,0,0,1,2) & 3 \\
3 & (0,0,0,1,1,1) & 3 \\
3 & (0,0,1,1,2,2) & 0 \\
\hline
4 & (0,0,2,2,1,3) & 0 
\end{tabular}
\qquad
\begin{tabular}{c||c|c}
 $n$ & $(k_1, \cdots, k_6)$ & $|\Tr \, g|$ \\
\hline
5 & (0,0,1,2,3,4) & 1 
\\ \hline
6 & (0,3,2,5,4,4) & 2 \\
6 & (3,3,0,0,2,4) & 1
\\ \hline
7 & (1,5,4,6,2,3) & 1 \\ 
\hline
8 & (0,4,2,6,1,3) & $\sqrt{2}$ 
\end{tabular}
\end{table}
where $g = \diag(\z_n^{k_1}, \cdots, \z_n^{k_6})$ and 
$|\Tr \, g| = |\z_n^{k_1} + \cdots + \z_n^{k_6}|$. 
Note that symplectic automorphisms of order $n \leq 8$ are conjugate to these normal 
forms as elements of $\PGL_6(\CC)$.  For example, let us consider $g = abab^{-1}$ 
with $a, b$ in Example \ref{exA7}. The order of $g$ is 2 in $\PGL_6(\CC)$ 
since we have $g^2 = \diag(\w^2, \cdots, \w^2)$. Eigenvalues of $g$ are
\[
\w, \ \w, \ \w, \ \w, \ -\w, \ -\w,
\]
and $g$ is conjugate to $\diag(1,1,1,1,-1,-1)$ in $\PGL_6(\CC)$.  
Also note that, although we can not define the trace of $g \in \PGL_6(\CC)$, the absolute value 
$|\Tr \, \tilde{g}|$ for a lift $\tilde{g} \in \GL_6(\CC)$ does not depend on a lifting if $\tilde{g}$ 
is of finite order. 
These values are usefull information to classify symplectic actions from character tables. 
\subsection{Invariants}
Let $G$ be one of Laza and Zheng's 34 groups. 
For our purpose, we need to consider 6-dimensional linear representations of $G^* \cong M.G$ such that 
$M$ acts as scalar matrices and the corresponding projective action of $G$ is faithful.  
In fact, we need only triple cover of $G$ as considered in \cite{YYZ23}. Let us see this in general 
situation.
In the following, $G$ denotes a finite group. 
\begin{Def}
For a homomorphism $\rho : G \rightarrow \GL_n(\CC)$ and $\chi \in \widehat{G}$, 
we say that a homogeneous polynomial $F(x) \in \CC[x_1, \dots, x_n]$ is a $\rho(G)$-semi-invariant 
with character $\chi$ if
\[
 F(\rho(g) x) = \chi(g) F(x) \quad \text{for} \quad {}^{\forall}g \in G.
\]
If $\chi$ is trivial, we say that $F(x)$ is a $\rho(G)$-invariant. \ We denote the vector space of 
$\rho(G)$-semi-invariants of degree $d$ with character $\chi$, 
by $V_d(\rho, \chi)$. If $\chi$ is trivial, it is denoted by $V_d(\rho)$. 
\end{Def}
\begin{rem} \label{Rem-chi}
{\rm (1)} \ In the above definition, it is obvious that $\chi$ factors through $\rho$, 
that is, there exists 
$\chi_{\rho} \in \widehat{\rho(G)}$ such that $\chi = \chi_{\rho} \circ \rho$. Moreover, 
$\chi_{\rho}$ is uniquely determined by $\chi$ and $\rho$.  
\\
{\rm (2)} \ 
For a Schur cover $G^* = M.G$ of $G$, we have $M \subset [G^*, G^*]$ and any homomorphism 
from $G^*$ to an abelian group factors through $G$. 
Therefore we have $\widehat{G^*} \cong \widehat{G}$.
\\
{\rm (3)} \ We have Molien's formula {\rm (see \cite{DK15})}
\begin{align*}
\frac{1}{|G|} \sum_{g \in G} \frac{1}{\det(\mathbf{1}_n - t \rho(g))} 
= \sum_{d=0}^{\infty} (\dim V_d(\rho))t^d,
\end{align*}
and we can compute $\dim V_d(\rho)$ using a function ``MolienSeries'' in GAP. 
\end{rem}
Now let $\rho : G \rightarrow \PGL_n(\CC)$ be a projective representation, and 
$\tilde{\rho} : G^* \rightarrow \GL_n(\CC)$ be a lifting of $\rho$, where $G^*$ is a Schur 
cover with the multiplier $M \subset G^*$. Assume that there is $\chi \in \widehat{G^*}$ 
such that $V_d(\tilde{\rho}, \chi) \ne 0$. 
\begin{lem} \label{SC-Lemma} 
In the above situation, we have
\\
{\rm (1)} \ $\tilde{\rho} (m)^d = \mathbf{1}_n$ for any $m \in M$.
\\
{\rm (2)} \ $\tilde{\rho}(M) \cong \C_d$ or $\tilde{\rho}(M) \cong 1$ if 
$d$ is prime.  
\end{lem}
\begin{proof}
{\rm (1)} \ By the assumption, $\tilde{\rho}(m)$ is a scalar matrix $\diag(c, \cdots, c)$. 
For a non-zero $F(x) \in V_d(\tilde{\rho}, \chi)$, we have
\[
F(\tilde{\rho}(m)x) =c^d F(x)
\]
and hence $\chi(m) = c^d$. Since $\mathrm{Im} \, \chi$ is abelian, we have 
\[
M \subset [G^*, G^*] \subset \mathrm{Ker} \, \chi. 
\]
Therefore we have $c^d = \chi(m) = 1$, that is, $c = (\z_d)^k$ for some $k \in \ZZ$. 
\\
{\rm (2)} follows from (1).   
\end{proof}
For example, the Schur multiplier of $\A_7$ is $\C_6$, and any projective representation 
$\rho : \A_7 \rightarrow \PGL_n(\CC)$ lifts to a linear representation 
$\tilde{\rho} : 6.\A_7 \rightarrow \GL_n(\CC)$. If $A_7$ acts on a cubic hypersurface in 
$\PP^{n-1}(\CC)$ via $\rho$, then $\tilde{\rho}$ factors through $3.\A_7 = 6.\A_7/\C_2$. 
\\ \\
We have the following simple facts that will be used to reduce the amount of calculation 
of invariants.
\begin{prop} \label{Prop-reduce}
Let $\rho : G \rightarrow \GL_n(\CC)$ be a representation. 
\\
{\rm (1)} \ For $\chi \in \widehat{G}$ and $\sigma \in \Aut(G)$, we have
\[
V_d(\rho, \chi) = V_d(\rho \circ \sigma, \, \chi \circ \sigma).
\] 
{\rm (2)} \ For $\chi, \, \theta \in \widehat{G}$, we have 
\[
V_d(\rho, \, \chi) = V_d(\rho \otimes \theta, \, \chi \otimes \theta^d).
\]
In particular, $\rho$ and $\rho \otimes \theta$ have same semi-invariants if 
$\theta^d = 1$.
\\ \\ 
{\rm (3)} \ Let us assume that $d$ and $|\widehat{G}|$ are co-prime to each other.
Then, for any $\chi \in \widehat{G}$, there exist $\theta \in \widehat{G}$ such that
\[
V_d(\rho, \, \chi) = V_d(\rho \otimes \theta).
\] 
Therefore, any semi-invariant is regarded as an invariant for a representation 
which is projectively equivalent to $\rho$.
\end{prop}
\begin{proof}
(1) and (2) follow from definitions. 
\\
(3) \ By the assumption, there is $e \in \NN$ such that 
\[
ed \equiv 1 \mod |\widehat{G}|.
\]
For $\theta = \chi^{-e}$, we have
\[
V_d(\rho, \, \chi) = V_d(\rho \otimes \theta, \, \chi \otimes \theta^d)
= V_d(\rho \otimes \theta, \, \chi^{1-ed}) = V_d(\rho \otimes \theta).
\]
\end{proof}
\subsection{Moduli spaces} \label{ModuliSpace}
Let us consider moduli spaces of hypersurfaces with action of finite group $G$ 
after \cite{YZ20}. For a representation $\rho : G \rightarrow \GL_n(\CC)$, 
we have the centralizer and the normalizer of $\rho(G)$ in $\GL_n(\CC)$:
\begin{align*}
\C_{\rho}(G) &= \{ h \in \GL_n(\CC) \ | \ hgh^{-1} = g \ \ \text{for} \ \ {}^{\forall}g \in \rho(G) \}, 
\\
\mathrm{N}_{\rho}(G) &= \{ h \in \GL_n(\CC) \ | \ h \rho(G) h^{-1} = \rho(G) \}. 
\end{align*}
Moreover we define the stabilizer of $V_d(\rho, \, \chi)$:
\begin{align*}
\mathrm{N}_{\rho, \chi}(G) &= \{ h \in \mathrm{N}_{\rho}(G) \ | \ 
\chi_{\rho}(h g h^{-1}) = \chi_{\rho}(g) \ \ \text{for} \ \ {}^{\forall}g \in \rho(G) \} 
\end{align*}
where $\chi_{\rho} \circ \rho = \chi$ (see Remark \ref{Rem-chi}).
The quotient group $\mathrm{N}_{\rho}(G) / \C_{\rho}(G)$ is isomorphic to a subgroup of 
$\Aut(\rho(G))$.  Hence the index $|\mathrm{N}_{\rho}(G) : \C_{\rho}(G)|$ is finite, and so is 
$|\mathrm{N}_{\rho, \chi}(G) : \C_{\rho}(G)|$. 
Let $V_d(\rho, \, \chi)^{sm}$ be the subset of semi-invariants which define smooth hypersurfaces of degree $d$. Then the quotient $V_d(\rho, \, \chi)^{sm} / \mathrm{N}_{\rho, \chi}(G)$ is the moduli 
space $\mathcal{M}_d(\rho, \chi)$ of smooth hypersurfaces with the action of $G$ with 
respect to $(\rho, \, \chi)$. 
\\ \indent
Now let $G$ be one of Laza and Zheng's 34 groups, $G^*$ be a Schur cover and 
$\rho : G^* \rightarrow \GL_6(\CC)$ be a lifting of a faithful projective action of $G$ on $\PP^5$.  
If $\rho(G^*)$ acts simplectically on a smooth cubic fourfold $F(x) = 0$ with $F \in V_3(\rho, \, \chi)$, 
we have
\begin{align} \label{Eq-dimension}
20 - r(G) = \dim \mathcal{M}_3(\rho, \chi) = \dim V_3(\rho, \, \chi) 
- \dim \mathrm{N}_{\rho, \chi}(G^*).
\end{align}
Note that we have $\dim \mathrm{N}_{\rho, \chi}(G^*) \geq 1$ since 
$\CC^{\times} \subset \mathrm{N}_{\rho, \chi}(G^*)$. Therefore we have the following criterion 
to exclude unnecessary cases. 
\begin{lem} \label{LM-dimension}
If we have
\[
\dim V_3(\rho, \, \chi) \leq 20 - r(G),
\]
then $F \in V_3(\rho, \, \chi)$ does not define a smooth cubic fourfold 
with a symplectic action of $\rho(G^*)$. 
\end{lem}
\begin{rem} \label{RM-moduli}
{\rm (1)} We have $\dim \mathrm{N}_{\rho, \chi}(G) = \dim \C_{\rho}(G)$, 
and $\dim \C_{\rho}(G)$ is easily computable. If we have a decomposition 
\[
\rho = m_1 \rho_1 + \cdots + m_k \rho_k
\]
into irreducible representations $\rho_i$, then we have
\[
\dim \C_{\rho}(G) = \dim (\GL_{m_1}(\CC) \times \cdots \times \GL_{m_k}(\CC)) 
= m_1^2 + \cdots + m_k^2
\]
by Schur's Lemma. 
\\
{\rm (2)} Obviously, we have
\[
\mathrm{N}_{\rho \otimes \chi}(G) = \mathrm{N}_{\rho}(G), \qquad
\mathrm{C}_{\rho \otimes \chi}(G) = \mathrm{C}_{\rho}(G)
\]
for $\chi \in \widehat{G}$.
\end{rem}
\subsection{$\A_6$-invariant cubic fourfolds} \label{SS-A6}
Let us consider $G = \A_6$ as an example. According to \cite{LZ22}, we have 1-dimensional family 
$\mathcal{M}(\A_6)$ of cubic fourfolds $X$ such that $\A_6 \subset \Aut^s(X)$, but it may have 
two or more irreducible components. 
\\ \indent
By Lemma (\ref{SC-Lemma}), an action of the Schur cover $6.\A_6$ on a cubic fourfold factors 
through $3.\A_6 = 6.\A_6/\C_2$ whose GAP-ID is (1080, 260). All we have to do is classify 
6-dimensional linear representations of $3.\A_6$ such that the center $\C_3$ acts as scalars. 
Since $\widehat{3.\A_6} = \widehat{\A_6}$ is trivial, we do not have a semi-invariant with 
a non-trivial character. Let us look at the character table of $3.\A_6$. 
\begin{align*}
\begin{array}{c|rrrrrrrrrrrrrrrrr} 
\hline
\rm class &\rm 1 &\rm 3A &\rm 3B &\rm 2A &\rm 6A  &\rm 6B &\rm 3C &\rm 3D &\rm 12A 
&\rm 12B &\rm 4A & \cdots 
\\
\rm size& 1 & 1 & 1 & 45 & 45 & 45 & 120 & 120 & 90 & 90 & 90 & \cdots
\\ \hline
\chi_1 & 1 & 1 & 1 & 1 & 1 & 1 & 1 & 1 & 1 & 1 & 1 & \cdots \\
\chi_2 & 3 & 3\bw & 3\w & -1 & -\bw & -\w & 0 & 0 & \bw & \w & 1 & \\
\chi_3 & 3 & 3\bw & 3\w & -1 & -\bw & -\w & 0 & 0 & \bw & \w & 1 & \\
\chi_4 & 3 & 3\w & 3\bw & -1 & -\w & -\bw & 0 & 0 & \w & \bw & 1 & \\
\chi_5 & 3 & 3\w & 3\bw & -1 & -\w & -\bw & 0 & 0 & \w & \bw & 1 & \\
\chi_6 & 5 & 5 & 5 & 1 & 1 & 1 & -1 & 2 & -1 & -1 & -1 & \\
\chi_7 & 5 & 5 & 5 & 1 & 1 & 1 & 2 & -1 & -1 & -1 & -1 & \\
\chi_8 & 6 & 6\bw & 6\w & 2 & 2\bw & 2\w & 0 & 0 & 0 & 0 & 0 & \\
\chi_9 & 6 & 6\w & 6\bw & 2 & 2\w & 2\bw & 0 & 0 & 0 & 0 & 0 & \\
\vdots & \vdots & \\
\hline
\end{array}
\end{align*}
\qquad \quad The following infomations are omitted: 
\\
\qquad \quad $\cdot$ irreducible representations whose dimensions are greater than 6,
\\
\qquad \quad $\cdot$ conjugacy classes for order $5$ and $15$. 
\\ \\
Classes 1, 3A and 3B form the center $\C_3$, and $\chi_i \ (i = 1, 6, 7)$ are irreducible 
characters with $\C_3 \subset \mathrm{Ker} \, \chi_i$, that is, characters of $\A_6$ itself 
rather than the triple cover $3.\A_6$. In this case, we have two (non-trivial) characters
$\chi_1 + \chi_6$ and $\chi_1 + \chi_7$ of degree 6. These characters are afforded by the following 
two representations. One is permutatons of projective coordinates $[x_1:x_2:x_3:x_4:x_5:x_6]$ 
denoted by $\rho$, and another is the composition of $\rho$ and an outer automorphism 
$\sigma$ of $\A_6$ (see Appendix A).
We have $V_3(\rho) = V_3(\rho \circ \sigma)$ (see Proposition \ref{Prop-reduce}) and this vector 
space is given by
\begin{align} \label{eqA6-1}
\mathrm{Span} \{ x_1^3 + \cdots + x_6^3, \quad (x_1^2 + \cdots + x_6^2)(x_1 + \cdots + x_6), 
\quad (x_1 + \cdots + x_6)^3 \}. 
\end{align}
This linear system contains the Fermat cubic, and a general member is smooth. 
 Since $\A_6$ is a non-commutative simple group, any action is symplectic  
(Although $\Sy_6$ acts on these cubics, the action of $(12) \in \Sy_6$ is not 
 symplectic since it is diagonalized into $\diag(1,1,1,1,1,-1)$, see Table \ref{TBL-Fu} ).
We have $\dim \C_{\rho}(3.\A_6) = 1^2 + 1^2 = 2$ (see Remark \ref{RM-moduli}) and 
\[
\dim V_3(\rho) - \dim \C_{\rho}(3.\A_6) = 1 = 20 -r(\A_6)
\]
as desired.
\\ \indent
Next we consider 6-dimensional representations such that a class 3A acts as a 
scalar $\w$ or $\w^2$. We have 
eight cases:
\[
2\chi_2, \quad \chi_2 + \chi_3, \quad 2\chi_3, \quad 2\chi_4, \quad \chi_4 + \chi_5, 
\quad 2\chi_5, \quad \chi_8, \quad \chi_9. 
\]
The image of class 4A in $\PGL_6(\CC)$ remains a class of order $4$, and (the absolute value of) 
the trace must be 0 if the action is symplectic (see Table \ref{TBL-Fu}). Since we have 
\[
\chi_i( 4\A) = \begin{cases}1 \quad (i = 2,3,4,5) \\ 0 \quad (i=8,9) \end{cases}, 
\] 
only $\chi_8$ and $\chi_9$ are admissible. From the character table, they are complex conjugate 
to each other. We give corresponding matrix representations $\rho$ and $\overline{\rho}$ as 
a subgroup of $G \cong 3.\A_7$ considered in Example \ref{exA7}.  
\\ \indent
For $a$ and $b$ in (\ref{ab in 3.A7}), using GAP, we can check that elements
\begin{align}
s = a^2(ba^2b^2a)^2 a, \quad t = b^2a^2ba 
\end{align}
generate a subgroup $H \cong 3.\A_6$. The $H$-invariant cubic forms are spanned by 
$f(x) $ in (\ref{3.A7-invariant}) and
\begin{multline} 
g(x) = 3 (x_1 x_2 x_3 + x_1 x_2 x_4 + x_1 x_2 x_5 + x_1 x_2 x_6 + x_1 x_3 x_4 + x_1 x_3 x_5 + 
   x_1 x_5 x_6 \\
   + x_2 x_3 x_4 + x_2 x_5 x_6 + x_3 x_5 x_6 + x_3 x_4 x_6 + x_3 x_4 x_5 + x_4 x_5 x_6) \\
   +4 (x_1 x_4 x_5 + x_1 x_3 x_6 + x_2 x_3 x_5 + x_2 x_4 x_6) \\
   +5 (x_1 x_4 x_6 + x_2 x_3 x_6 + x_2 x_4 x_5) + 8 x_1 x_3 x_5.
\end{multline} 
This linear system contains a smooth cubic fourfolds $f(x)=0$, and a general member 
is smooth. Moreover we have $H = \overline{H}$ since 
\[
\overline{s} = t^2 s t^{-1} (st)^2 t, \qquad 
\overline{t} = s t (st^{-1})^2 s t s t^{-1} s.
\]
Namely we have
\begin{align} \label{eqA6-2}
V(\rho) = V(\overline{\rho}) = \CC f(x) \oplus \CC g(x).
\end{align}
From the above, we see that $\mathcal{M}(A_6)$ has two irreducible components.
\subsection{Proof of Theorem 1.1} \label{proof}
As is the case of $\A_6$, we determine all admisible projective representations 
(excluding $\rho$ for which $V_3(\rho)$ is same with invariants already obtained) for
 34 groups $G$ in the following sections.  The result is listed in Table \ref{TBL-Schur} in 
 Appendix B. To complete the proof of Theorem 1.1,  we need to show that two 
 families for groups
 \begin{align*}
3, \ \A_{3,3}, \ \D_{12}, \ \A_4, \ \A_{4,3}, \ \A_5, \ 3^2.4, \ \A_6, \ \Sy_5, 
\ \A_7, \ \M_{10}
\end{align*}
 are distinct (see \cite{LZ22} for $3, \ \D_{12}, \ \A_7, \ \M_{10}$). This is shown as follows. 
 Note that if $\psi : X_1 \rightarrow X_2$ is an isomorphism of smooth cubic fourfolds, 
 then it induces a group isomorphism $\phi : \Aut^s(X_1) \rightarrow \Aut^s(X_2), \ 
 \phi(g) = \psi \circ g \circ \psi^{-1}$ with a commutative diagram
 \[
 \begin{CD} 
 \Aut^s(X_1) \times X_1 @>>> X_1 \\
 @V{(\phi, \psi)}VV @VV{\psi}V\\
 \Aut^s(X_2) \times X_2 @>>> X_2
 \end{CD}
 \] 
 where horizontal arrows are group actions. Therefore if two family $V_3(\rho_1, \chi_1)$ 
 and $V_3(\rho_2, \chi_2)$ for $G$ listed in the above are same, then there exist 
 $\phi \in \Aut(G)$ such that $\rho_1$ is projectively equivalent to $\rho_2 \circ \phi$. 
 However, we see that this is impossible by considering actions of Schur multiplier or 
comparing dimensions of irreducible components of two representations. 
For example, we have two families (\ref{eqA6-1}) and (\ref{eqA6-2}) for $G = \A_6$. 
We can distinguish them from each other since the Schur multiplier $\C_3$ acts trvially 
for (\ref{eqA6-1}) and non-trivially for (\ref{eqA6-2}). We see also that irreducible decomposition 
of the action of $\A_6$ consists of a 1-dimensional representation and  a 5-dimensional 
representatoin for (\ref{eqA6-1}) but an irreducible 6-dimensional representation for (\ref{eqA6-2}). 

%% file: SmallGroup.tex
\section{Groups of small $r(G)$}
\subsection{Groups of $r(G) \leq 14$}
Let $G$ be one of small groups
\[
\C_2, \quad \C_2 \times \C_2, \quad \C_3, \quad \C_4, \quad \Sy_3.
\]
Cubic fourfolds $X$ such that $G \subset \Aut^s(X)$ are determined in \cite{Fu16} 
and \cite{LZ22} as follows.  
\\ \\
(i) \ $G = \C_2$. We may assume that $G$ is generated by $\diag(1,1,1,1,-1,-1)$, 
and invariant cubics are 
\begin{multline} \label{eqC2}
F_1(x_1, x_2, x_3, x_4) + x_5^2 L_1(x_1, x_2, x_3, x_4) \\
+ x_5 x_6 L_2(x_1, x_2, x_3, x_4) + x_6^2 L_3(x_1, x_2, x_3, x_4) = 0.
\end{multline}
(ii) \ $G = \C_2 \times \C_2$. We may assume that $G$ is generated by
$\diag(1,1,1,1,-1,-1)$ and $\diag(1,1,1,-1,-1,1)$, and invariant cubics are
\begin{multline} \label{eqC2xC2}
F_1(x_1, x_2, x_3) + x_4^2 L_1(x_1, x_2, x_3) + x_5^2 L_2(x_1, x_2, x_3) \\
+x_6^2 L_3(x_1, x_2, x_3) + x_4 x_5 x_6 = 0
\end{multline}
(iii) \ $G = \C_3$. There are two families. One is 
\begin{align} \label{eqC3-1}
F_1(x_1, x_2, x_3, x_4) + x_5^3 + x_6^3 + x_5 x_6 L_1(x_1, x_2, x_3, x_4) = 0
\end{align}
and $G$ is generated by $\diag(1,1,1,1,\w,\w^2)$. Another is 
\begin{align} \label{eqC3-2}
F_1(x_1, x_2) + F_2(x_3, x_4) + F_3(x_5, x_6) 
+ \sum_{i=1,2 ; j=3,4 ; k=5,6} a_{ijk} x_i x_j x_k = 0
\end{align}
and $G$ is generated by $\diag(1,1,\w,\w,\w^2,\w^2)$.
\\
(iv) \ $G = \C_4$. We may assume that $G$ is generated by $\diag(1,1,-1,-1,i,-i)$
and invariant cubics are
\begin{multline} \label{eqC4}
\mathrm{Span} \{ x_1 N_1(x_3, x_4),\ x_2 N_2(x_3, x_4),\ F_1(x_1, x_2), \\
x_5 x_6 L_1(x_1, x_2),\ x_5^2 L_2(x_3, x_4),\ x_6^2 L_3(x_3, x_4) \}.
\end{multline}
(v) \ $G = \Sy_3$. We have the following normal forms. The symmetric group 
$\Sy_3$ acts on projective coordinates
\[
[x_1 : x_2 : x_3 : y_1 : y_2 : y_3]
\] 
by permuting $x_i$ and $y_i$ simultaneously. Invariant cubics are 
\begin{multline} \label{eqS3}
\mathrm{Span} \{ s_1(x)^3, \ s_1(x) s_2(x), \ s_3(x), \ s_1(y)^3, \ s_1(y)s_2(y), \ s_3(y), \\ 
s_1(x)s_1(y)^2, \ s_1(x)s_2(y), \ s_1(x)^2s_1(y), \ s_2(x)s_1(y), \\
\sum_{i=1}^3 x_iy_i^2, \ \sum_{i=1}^3 x_i^2y_i, \ 
\sum_{\sigma \in \A_3} x_{\sigma(1)} y_{\sigma(2)} y_{\sigma(3)},
\ \sum_{\sigma \in \A_3} x_{\sigma(1)} x_{\sigma(2)} y_{\sigma(3)} \}
\end{multline} 
where $s_i(x) = x_1^i + x_2^i + x_3^i$ and $\A_3$ is the alternating group of degree 3. 
\subsection{The dihedral group $\D_8$}
The Schur multiplier of $\D_8$ is $\C_2$, and we do not need a Schur cover 
(Lemma \ref{SC-Lemma} ). As an abstract group, $\D_8$ is 
\[
\D_8 = \left< a, b \ | \ a^4 = b^2 = 1, \ bab = a^{-1} \right>,
\]
and the charactet table is as follows.
\[
\begin{array}{c|rrrrr}
\hline
 \rm class&1& \rm2A & \rm2B\quad & \rm2C\quad & \rm4A\quad \cr
  &\rm1& a^2 &\{ b, ba^2\}&\{ ba, ba^3\}&\{a, a^3 \}\cr
\hline
  \chi_{1}&1&1&1\quad&1\quad&1\quad \cr
  \chi_{2}&1&1&-1\quad&-1\quad&1\quad \cr
  \chi_{3}&1&1&1\quad&-1\quad&-1\quad \cr
  \chi_{4}&1&1&-1\quad&1\quad&-1\quad \cr
  \chi_{5}&2&-2&0\quad&0\quad&0\quad \cr
\hline
\end{array}.
\]
We have $\widehat{\D_8} = \{ \chi_1, \chi_2, \chi_3, \chi_4\} \cong \C_2 \times \C_2$, 
and we do not have to consider semi-invariants by Proposition \ref{Prop-reduce} (3). 
If $\rho : \D_8 \rightarrow \GL_6(\CC)$ satisfies the symplectic conditions, we have
\[
|\Tr \rho(g)| = \begin{cases} 2 \quad (\text{the class of} \ g \ \text{is} \ 
2\mathrm{A},  2\mathrm{B}, 2\mathrm{C}) \\ 
0 \quad (\text{the class of} \ g \ \text{is} \ 4\mathrm{A}) \end{cases},
\]
and these conditions are satisfied only by
\begin{align*}
\theta_1 = 2\chi_1 + \chi_3 + \chi_4 + \chi_5, \qquad
\theta_2 = 2\chi_2 + \chi_3 + \chi_4 + \chi_5, \\
\theta_3 = \chi_1 + \chi_2 + 2\chi_3 + \chi_5, \qquad
\theta_4 = \chi_1 + \chi_2 + 2\chi_4 + \chi_5.
\end{align*}
The outer automorphism group $\Out(\D_8) = \C_2$ is generated by 
\[
\sigma : \D_8 \longrightarrow \D_8, \qquad \sigma(a) = a^3, \quad \sigma(b) = ab.  
\]
It interchanges classes 2B and 2C, and hence $\chi_3$ and $\chi_4$. Therefore invariants for  
$\theta_3$ and $\theta_4$ are same. Let us construct representations $\rho_i$ 
affording $\theta_i \ (i=1,2,3)$.
\\ \indent
We may assume that $a = \diag(1,1,-1,-1, i, -i)$, and then $\rho_i$ is given by
$a$ and the following $b_i$:  
\begin{align*}
b_1 &= \diag(1,1,1,-1) \oplus \begin{bmatrix} 0 & 1 \\ 1 & 0\end{bmatrix}, \quad
b_2 = \diag(-1,-1,1,-1) \oplus \begin{bmatrix} 0 & 1 \\ 1 & 0\end{bmatrix}, \\
b_3 &= \diag(1,-1,1,1) \oplus \begin{bmatrix} 0 & 1 \\ 1 & 0\end{bmatrix},
\end{align*}
which satisfy symplectic conditions. Now $V_3(\rho_i)$ is given as a subspace of 
$b_i$-invariants in (\ref{eqC4}). We have
\begin{multline} \label{eqD8}
V_3(\rho_1) = \mathrm{Span}\{ x_1 x_3^2, \ x_1 x_4^2, \ x_2 x_3^2, \ x_2 x_4^2, \
F_1(x_1, x_2),\ x_5 x_6 L_1(x_1, x_2), \\
(x_5^2 + x_6^2) x_3,\ (x_5^2 - x_6^2) x_4 \},
\end{multline}
\begin{align}
V_3(\rho_2) = \mathrm{Span}\{ x_1 x_3 x_4, \  x_2 x_3 x_4, \ (x_5^2 - x_6^2)x_4,\
(x_5^2 + x_6^2) x_3 \},
\end{align}
\begin{align}
V_3(\rho_3) = \mathrm{Span}\{ x_1 N_1(x_3, x_4), \ x_1^3, \ x_1 x_2^2, \ 
x_5 x_6 x_1, \ (x_5^2 + x_6^2) L_2(x_3, x_4) \}.
\end{align}
A general member is a smooth cubic fourfold only for $V_3(\rho_1)$. 
A cubic fourfold defined by $F \in V_3(\phi_2)$ contains $\PP^3 = \{x_3 = x_4 =0\}$ 
and singular. Similarly, $F \in V_3(\phi_3)$ vanishes on 
$\PP^3 = \{x_1 = x_5 + ix_6 =0\}$. We have 
\[
\dim \C_{\rho_1}(\D_8) = 2^2 + 1^2 + 1^2 + 1^2 = 7
\]
and
\[
\dim V_3(\rho_1) - \dim \C_{\rho_1}(\D_8) = 5 = 20 - r(\D_8).
\]

%% file: r16.tex
\section{Groups of $r(G) = 16$}
We have the following four groups
\[
 \A_{3,3}, \quad \D_{12}, \quad \A_4, \quad \D_{10}.
\]
\subsection{The dihedral groups $\D_{10}$ and $\D_{12}$} 
For $G = \D_{10}, \ \D_{12}$, cubic fourfolds such that $G \subset \Aut^s(X)$ 
are determined in \cite{LZ22} as follows.
\\
(i) \ $G = \D_{10}$. For an appropriate choice of coordinates,
\begin{multline} \label{eqD10}
F_1(x_1, x_2) + x_3 x _6 L_1(x_1, x_2) + x_4 x_5 L_2(x_1, x_2) \\
+ x_3^2 x_5 +x_3 x_4^2 + x_4 x_6^2 + x_5^2 x_6 = 0.
\end{multline}
An order 5 element in $G$ is $g = \diag(1, 1, \z_5, \z_5^2, \z_5^3, \z_5^4)$. 
Moreover, any smooth cubic fourfolds with a symplectic automorphism of order 5 have 
this form, and a generic such cubic fourfold has symplectic automorphism group $\D_{10}$ 
which is generated by $g$ and a permutation coordinates by $(36)(45)$.
\\ \\
(ii) \ $G = \D_{12}$. The defining equations of the corresponding cubic fourfolds either
belong to
\begin{multline} \label{eqD12-1}
\mathrm{Span} \{ x_1^2 x_3, x_1^2 x_4, x_1 x_2 x_3, x_1 x_2 x_4, 
x_2^2 x_3, x_2^2 x_4, x_3^3, \\
x_3^2 x_4, x_3 x_4^2, x_3 x_5 x_6, x_4^3, x_4 x_5 x_6, x_5^3, x_6^3 \},
\end{multline}
while an order 6 element of $G$ is $\diag(-1,-1,1,1,\w, \w^2)$ or belong to
\begin{multline} \label{eqD12-2}
\mathrm{Span} \{x_1^3, x_1 x_2^2, x_1 x_3 x_5, x_1x_3x_6, x_2x_4x_5, 
x_2x_4x_6, x_3^3 , x_3x_4^2, x_5^3 , x_5^2 x_6, x_5x_6^2, x_6^3 \},
\end{multline}
while an order 6 element of $G$ is $\diag(1,-1,\w,-\w,\w^2,\w^2)$. 
Moreover, a generic cubic fourfold admitting such an order 6 automorphism has symplectic 
automorphism group $\D_{12}$.
\\ \indent
The action of $\D_{12}$ is given by the following way. A general member of (\ref{eqD12-1}) is 
normalized into
\begin{align}
a_1(x_1^2 x_3 + x_2^2 x_3) + a_2(x_1^2 x_4 + x_2^2x_4) + a_3(x_5^6 + x_6^6) + \text{other terms}
\end{align}
by diagonal actions of $(\CC^{\times})^6$. Then $\D_{12}$ is generated by 
$\diag(-1,-1,1,1,\w, \w^2)$ and a permutation of coordinates by $(12)(56)$.
\\ \indent
A general member of (\ref{eqD12-2}) is normalized into
\begin{align}
a_1(x_1^3 + x_3^3) + a_2(x_1 x_2^2 + x_3 x_4^2) + \text{other terms}
\end{align}
by diagonal actions of $(\CC^{\times})^6$. Then $\D_{12}$ is generated by 
$\diag(1,-1,\w,-\w,\w^2,\w^2)$ and a permutation of coordinates by $(13)(24)$.
\subsection{The groups $\A_{3,3} =  (\Sy_3 \times \Sy_3) \cap \A_6$} 
\label{SS-A33}
The Schur multiplier of $G = \A_{3,3}$ is $\C_3$, and we consider a Schur cover 
$G^* = 3.\A_{3,3}$ whose GAP-ID is (54,8). 
Firstly, we consider linear actions of $A_{3,3}$ itself. 
The characters teble of $A_{3,3}$ is as follows.
\begin{align*}
\begin{array}{c|rrrrrr}
  \rm class&1& \rm2A & \rm3A & \rm3B & \rm3C & \rm3D \cr
  \rm size&1&9&2&2&2&2\cr
\hline
  \chi_{1}&1&1&1&1&1&1\cr
  \chi_{2}&1&-1&1&1&1&1\cr
  \chi_{3}&2&0&2&-1&-1&-1\cr
  \chi_{4}&2&0&-1&2&-1&-1\cr
  \chi_{5}&2&0&-1&-1&2&-1\cr
  \chi_{6}&2&0&-1&-1&-1&2
\end{array}
\quad 
\begin{array}{l}
\rm 2A : (12)(45) \\
\rm 3A : (123) \\
\rm 3B : (456) \\
\rm 3C : (123)(456) \\
\rm 3D : (123)(465) \\
\end{array}
\end{align*}
Since we have $\widehat{A_{3,3}} \cong C_2$, we do not have to consider semi-invariants. 
If $\rho : \A_{3,3} \rightarrow \GL_6(\CC)$ is faithful and satisfies the symplectic conditions, 
we have
\[
|\Tr_{\rho}(g)| = \begin{cases} 2 \quad (\text{the order of $g$ is 2}) \\ 
0, \, 3 \quad (\text{the order of $g$ is 3}) \end{cases},
\]
and these conditions are satisfied only by
\[
\theta_{ij} = 2\chi_1 + \chi_i + \chi_j, \quad \theta_{ij} \otimes \chi_2  \qquad (3 \leq i < j \leq 6). 
\]
The outer automorphism group $\Out(A_{3,3})$ is isomorphic to $\Sy_4$, and it permutes 
classes $\rm3A, \ 3B, \ 3C$ and $\rm3D$. We consider only $\theta_{34}$ and 
$\theta_{34} \otimes \chi_2$ since $\Out(\A_{3,3})$ permutes $\theta_{ij}$ 
(see Proposition \ref{Prop-reduce} (1)). 
\begin{rem} We see that $\Out(A_{3,3})$ is generated by conjugations by 
$(12),(13)(24)(36)$, and an automorphism $\sigma$ in Appendix A.
\end{rem}
A representation $\rho$ affording 
$\theta_{34}$ is given by permutations of coordinates $[x_1 : x_2 : x_3 : y_1 : y_2 : y_3]$, 
and $V_3(\rho)$ is 
\begin{multline} \label{eqA33-1}
\mathrm{Span} \{ s_1(x)^3, \ s_2(x)s_1(x), \ s_3(x), \ s_1(y)^3, \ s_2(y)s_1(y), \ s_3(y), \\
s_1(x)s_2(y), \ s_1(y)s_2(x), \ s_1(x)s_1(y)^2, \ s_1(x)^2 s_1(y) \},
\end{multline}
which is a subspace of (\ref{eqS3}). Although these are stable under permutatoins by 
$\Sy_{3,3}$, the action of $(12) \in \Sy_{3,3}$ is non-symplectic. 
We have
\[
 \dim \C_{\rho}(\A_{3,3}) = 2^2 + 1^2 + 1^2 = 6
\]
and
\[
\dim V_3(\rho) - \dim \C_{\rho}(\A_{3,3}) = 4 = 20 - r(\A_{3,3}).
\]
By Proposition \ref{Prop-reduce} (2), we have $V_3(\rho \otimes \chi_2) = V_3(\rho, \chi_2)$ 
and it is
\begin{align} \label{eqA33-singular}
 \mathrm{Span} \{ (x_1 - x_2)(x_1 - x_3)(x_2 - x_3), \
  (y_1 - y_2)(y_1 - y_3)(y_2 - y_3) \}.
\end{align}
This is not the case. 
\\ \indent
Next, we consider $G^* = 3.A_{3,3}$. Since we have $\widehat{G^*} \cong \widehat{G} \cong \C_2$, 
 we do not have to consider semi-invariants. Irreducible representaions of $G^*$ 
with a non-trivial action of the Schur multiplier $\C_3$ are given as follows.
\begin{align*}
\begin{array}{c|rrrrrrrrrr}
\hline
  \rm class&\rm1&\rm2A&\rm3A&\rm3B&\rm3C&\rm3D&\rm3E&\rm3F&\rm6A&\rm6B\cr
  \rm size&1&9&1&1&6&6&6&6&9&9\cr
\hline
  \chi_{7}&3&1& 3\bw & 3\w &0&0&0&0&\w&\bw \cr
  \chi_{8}&3&1& 3\w & 3\bw &0&0&0&0&\bw&\w \cr
  \chi_{9}&3&-1& 3\w& 3\bw &0&0&0&0&-\bw&-\w \cr
  \chi_{10}&3&-1& 3\bw & 3\w &0&0&0&0&-\w&-\bw \cr
\hline
\end{array}
\end{align*}
Classes 1, 3A and 3B form the Schur multiplier, and they act as scalar matrices for
\[
 2\chi_i \ (i=7,8,9,10), \quad \chi_7 + \chi_{10}, \quad \chi_8 + \chi_9. 
\]
Since we have
\[
 (\chi_7 + \chi_{10})(2\A) = (\chi_8 + \chi_9)(2\A) = 0, 
\]
$\chi_7 + \chi_{10}$ and $\chi_8 + \chi_9$ are not suitable. Let $\rho_i$ be 
a representation affording $\chi_i$.  We have Molien series
\begin{align*}
 \mu(2\rho_i, \, t) = \sum_{d=0}^{\infty} (\dim V_d(2 \rho_i))t^d = 1 + 8 t^3 + \cdots \qquad (i=7, 8), \\
 \mu(2\rho_i, \, t) = \sum_{d=0}^{\infty} (\dim V_d(2 \rho_i))t^d = 1 + 42 t^6 + \cdots \qquad (i=9, 10), 
\end{align*}
and there is no invariant cubic for $2\rho_9$ and $2\rho_{10}$. Moreover, we have 
$\chi_8 = \overline{\chi_7}$,  and hence $\rho_8$ is the complex conjugate of $\rho_7$. 
\\ \indent
Let us consider matrices
\begin{align} \label{abc-A33}
a = \begin{bmatrix}
1 & 0 & 0 \\ 0 & 0 & 1 \\ 0 & 1 & 0  
\end{bmatrix},\quad 
b = \begin{bmatrix}
0 & \w & 0 \\ \w^2 & 0 & 0 \\ 0 & 0 & 1
\end{bmatrix},\quad
c = \begin{bmatrix}
1 & 0 & 0 \\ 0 & \w & 0 \\ 0 & 0 & \w^2
\end{bmatrix}
\end{align}
with eigenvalues $(1,1,-1),\ (1,1,-1),\ (1, \w, \w^2)$. These matrices generate 
a group $H$ with GAP-ID (54,8), and this is one of $\rho_7$ and $\rho_8$. 
Note that we have 
\begin{align} \label{abc-conj}
\overline{a} = a, \quad \overline{b} = c^2bc,\quad \overline{c} = aca
\end{align}
and hence $\overline{H} = H$, that is, $V_3(2\rho_7) = V_3(2\rho_8)$. 
From eigenvalues, we see that matrices 
\begin{align} \label{abc2-A33}
\tilde{a} = \begin{bmatrix} a & {\bf 0} \\ {\bf 0} & a \end{bmatrix},\quad 
\tilde{b} = \begin{bmatrix} b & {\bf 0} \\ {\bf 0} & b \end{bmatrix},\quad 
\tilde{c} = \begin{bmatrix} c & {\bf 0} \\ {\bf 0} & c \end{bmatrix},\quad 
\end{align}
satisfy symplectic conditions, and $V_3(2\rho_7)$ with projective coordinates 
$[x_1:x_2:x_3:y_1:y_2:y_3]$ is given by
\begin{align} \label{eqA33-2}
{\rm Span} \{ s_3(x), \ x_1 x_2 x_3, \ s_3(y), \ y_1 y_2 y_3, \ I_1, \ I_2, \ I_3, \ I_4 \}
\end{align}
where 
\begin{align*}
I_1 &= x_1 y_1^2 + x_2 y_2^2 + x_3 y_3^2, \qquad
I_2 = x_1 y_2 y_3 + y_1 x_2 y_3 + y_1 y_2 x_3 \\
I_3 &= x_1^2 y_1 + x_2^2 y_2 + x_3^2 y_3, \qquad
I_4 = y_1 x_2 x_3 + x_1 y_2 x_3 + x_1 x_2 y_3. 
\end{align*}
We have $\dim \C_{\rho_7}(3.\A_{3,3}) = 2^2 = 4$ and
\[
\dim V_3(\rho_7) - \dim \C_{\rho_7}(3.\A_{3,3}) = 4 = 20 - r(\A_{3,3}).
\]
\begin{rem}  
We can also consider $3.A_{3,3}$ as a subgroup of $H \cong 3.\A_6$ 
in section {\rm \ref{SS-A6}}. Here we chose the above matrices 
(which is not in $H$) to obtain invariants of simple forms.  
\end{rem}
\subsection{The alternating group $\A_4$} 
The Schur multiplier of $\A_4$ is $\C_2$, and we do not need a Schur cover. 
The character table of $\A_4$ is as follows:
\\
\begin{align} \label{CT-A4}
\begin{array}{c|cccc}
  \rm class& 1 & (12)(34) & (123) & (124) \cr
  \rm size&1&3&4&4\cr
\hline
  \chi_{1}&1&1&1&1\cr
  \chi_{2}&1&1&\w &\bw \cr
  \chi_{3}&1&1&\bw&\w \cr
  \chi_{4}&3&-1&0&0\cr
\end{array}
\end{align}
We see that only the following characters satisfy the condition on absolute values:
\begin{align} \label{REP-A4}
\theta_i = 3 \chi_i + \chi_4 \ (i = 1,2,3), \quad \theta_4 = 2 \chi_4, \quad 
\theta_5 = \chi_1 + \chi_2 + \chi_3 + \chi_4
\end{align}
(Note that we have $\theta_i = \theta_1 \otimes \chi_i$ for $i = 1, 2, 3$). 
Third symmetric powers of these characters are
\begin{align*}
\sym^3 \theta_i &= 14 \chi_1 + 3 \chi_2 + 3 \chi_3 + 12 \chi_4 \ (i =1,2,3), \\
\sym^3 \theta_4 &= 4 \chi_1 + 2 \chi_2 + 2 \chi_3 + 16 \chi_4, \\
\sym^3 \theta_5 &= 8 \chi_1 + 6 \chi_2 + 6 \chi_3 + 12 \chi_4.  
\end{align*}
Let $\rho_i$ be a representation affording $\theta_i$. Then we have
\begin{align*}
\dim V_3(\rho_i) = 14 \ (i=1,2,3), \quad 
\dim V_3(\rho_5) = 8, \qquad \dim V_3(\rho_5, \chi_i) = 6 \ (i=2,3)
\end{align*}
and we can expect smooth cubics with symplectic action of $\A_4$  
only for these spaces by Lemma \ref{LM-dimension}.  
By Proposition \ref{Prop-reduce} (2), we have
\[
V_3(\rho_1) = V_3(\rho_1 \otimes \chi_i, \, \chi_i^3) 
= V_3(\rho_1 \otimes \chi_i) = V_3(\rho_i) 
\]
since $\chi_i^3 = 1$ for $i = 1,2,3$. 
\\ \indent
Taking projective coordinates $[x_1 : x_2 : x_3 : x_4 : y_1 : y_2]$, the representation $\rho_1$ 
is given by permutations of $x_i$, and it satisfies the symplectic condition. The linear 
space $V_3(\rho_1)$ is 
\begin{multline} \label{eqA4-1}
\mathrm{Span}\{ s_1(x)s_2(x), \ s_2(x) y_1, \ s_2(x) y_2, \ s_3(x), \
 \text{cubic forms in} \ s_1(x), y_1, y_2 \}
\end{multline}
where $s_i(x) = x_1^i + x_2^i + x_3^i + x_4^i$. A general member of $V_3(\rho_1)$ is smooth 
since it contains the Fermat cubic. 
We have $ \dim \rho_1(\A_4) = 3^2 + 1^2 = 10$ and 
\[
\dim V_3(\phi_1) - \dim \C_{\rho_1}(\A_4) = 4 = 20 - r(\A_4).
\]
The representation $\rho_5$ is given by
\begin{align} \label{ACT-3.A4}
\begin{split}
(123)[x_1 : x_2 : x_3 : x_4 : y_1 : y_2] = [x_2 : x_3 : x_1 : x_4 :\w y_1 :\w^2 y_2],  \\
(12)(34)[x_1 : x_2 : x_3 : x_4 : y_1 : y_2] = [x_2 : x_1 : x_4 : x_3 : y_1 : y_2], 
\end{split}
\end{align}
and it satisfies symplectic conditions. The linear space $V_3(\rho_5)$ is spanned by
\begin{align} \label{eqA4-2}
s_1(x)^3, \quad s_1(x)s_2(x), \quad s_3(x), \quad s_1(x)y_1y_2, \quad
y_1^3, \quad y_2^3, \quad f(x)y_2, \quad g(x)y_1,
\end{align}
where
\begin{align*}
f(x) &= (x_1-x_2)(x_3-x_4) -\w(x_1-x_3)(x_2-x_4) +\w^2(x_1-x_4)(x_2-x_3), \\
g(x) &= (x_1-x_2)(x_3-x_4) -\w^2(x_1-x_3)(x_2-x_4) +\w(x_1-x_4)(x_2-x_3).
\end{align*}
A general member of $V_3(\rho_5)$ is smooth since it contains the Fermat cubic.  
We have 
\[
\dim \C_{\rho_5}(\A_4) = 1^2 + 1^2 + 1^2 + 1^2 = 4 \]
 and 
\[
\dim V_3(\rho_5) - \dim \C_{\rho_5}(\A_4) = 4 = 20 - r(\A_4).
\]
Since we have
\[
\dim V_3(\rho_5, \ \chi_i) - \dim \C_{\rho_5}(\A_4) = 6-4 = 2 \quad (i = 2,3),
\]
a general member of $V_3(\rho_5, \chi_i)$ is not smooth for $i =2,3$.

%% file: r17.tex
\section{Groups of $r(G) = 17$}
We have two groups $\Sy_4$ and $\Q_8$.
\subsection{The quaternion group $\Q_8$} \label{SS-Q8}
The Schur multiplier of $\Q_8$ is 1, and hence we have $\Q_8^* = \Q_8$.  We express $\Q_8$ as
\[
\Q_8 = \{ \pm1, \ \pm I, \ \pm J, \ \pm K\}
\]
with usual rules $I^2 = J^2 = K^2 = -1$, \ $IJ = -JI = K$. 
The character table of $\Q_8$ is the following:
\begin{align*}
\begin{array}{c|rrrrr}
  \rm class&\rm1& -1 & \pm I & \pm J & \pm K\cr
\hline
  \chi_{1}&1&1&1&1&1\cr
  \chi_{2}&1&1&-1&1&-1\cr
  \chi_{3}&1&1&1&-1&-1\cr
  \chi_{4}&1&1&-1&-1&1\cr
  \chi_{5}&2&-2&0&0&0\cr
\end{array}.
\end{align*}
Since we have $\widehat{\Q_8} = \C_2 \times \C_2$, we do not have to consider semi-invariants. 
If $\rho : \Q_8 \rightarrow \GL_6(\CC)$ satisfies the symplectic conditions, we have
\[
|\Tr_{\rho}(g)| = \begin{cases} 2 \quad (g = -1) \\ 
0 \quad (g = I, J, K) \end{cases},
\]
and these conditions are satisfied only by
\[
\theta = \chi_1 + \chi_2 + \chi_3 + \chi_4 + \chi_5. 
\]
A representation $\rho$ affording $\theta$ is given by
\begin{align*}
\rho(I) = \diag(1,-1,1,-1) \oplus \begin{bmatrix} 0 & 1 \\ -1 & 0\end{bmatrix}, \quad
\rho(J) = \diag(1,1,-1,-1) \oplus \begin{bmatrix} i & 0 \\ 0 & -i\end{bmatrix}.
\end{align*}
Note that $\rho(J)$ is the same with the action of $\C_4$ for (\ref{eqC4}). 
The linear space $V_3(\rho)$ is given as $\rho(I)$-invariants in (\ref{eqC4}), and 
spanned by
\begin{align} \label{eqQ8}
 x_1^3, \ x_1x_2^ 2, \ x_1x_3^2, \ x_1 x_4^2, \ x_2 x_3 x_4, \ x_2 x_5 x_6, \ 
x_3(x_5^2 + x_6^2), \ x_4(x_5^2 - x_6^2).
\end{align}
We have 
\[
\dim \C_{\rho}(\Q_8) = 1^2 + 1^2 + 1^2 + 1^2 + 1^2 = 5
\]
and
\[
\dim V_3(\rho) - \dim \C_{\rho}(\Q_8) = 3 = 20 - r(\Q_8).
\]
\subsection{The symmetric group $\Sy_4$}
The Schur multiplier of $\Sy_4$ is $\C_2$, and we do not have to consider a Schur cover. 
The character tables of $\Sy_4$ is as follows:
\[
\begin{array}{c|rrrrr}
  \rm class & 1 & (12)(34) & (12) & (123) & (1234) \cr
  \rm size&1&3 \quad & 6 \ \ & 8 \quad & 6 \quad \cr
\hline
  \chi_{1}&1&1 \quad &1 \ & 1 \quad & 1 \quad \cr
  \chi_{2}&1&1 \quad &-1 \ & 1 \quad & -1 \quad \cr
  \chi_{3}&2&2 \quad &0 \ & -1 \quad & 0 \quad \cr
  \chi_{4}&3&-1 \quad &-1 \ & 0 \quad & 1 \quad \cr
  \chi_{5}&3&-1 \quad &1 \ & 0 \quad & -1 \quad \cr
\end{array}
\]
Since we have $\widehat{\Sy_4} = \C_2$, we do not have to consider semi-invariants. 
If $\rho : \Sy_4 \rightarrow \GL_6(\CC)$ satisfies the symplectic condition, we have
\[
|\Tr_{\rho}(12)(34)| = |\Tr_{\rho}(12)| =2, \quad |\Tr_{\rho}(123)| = 0,3, \quad 
|\Tr_{\rho}(1234)| = 0.
\]
These conditions are satisfied only for
\begin{align*}
\theta_1 = 2\chi_1 + \chi_2 + \chi_5, \quad \theta_2 =\theta_1 \otimes \chi_2, \quad
\theta_3 = \chi_1 + \chi_3 + \chi_5, \quad 
\theta_4 = \theta_3 \otimes \chi_2,
\end{align*}
and their third symmetric powers are
\begin{align*}
\sym^3 \theta_1 = 9 \chi_1 + 5 \chi_2 + 3 \chi_3 + 4 \chi_4 + 8 \chi_5, \quad
\sym^3 \theta_2 = \chi_2 \otimes \sym^3 \theta_1, \\
\sym^3 \theta_3 = 6 \chi_1 + 2 \chi_2 + 6 \chi_3 + 4 \chi_4 + 8 \chi_5, \quad
\sym^3 \theta_4 = \chi_2 \otimes \sym^3 \theta_3.
\end{align*}
Let $\rho_i : \Sy_4 \rightarrow \GL_6(\CC)$ be a representation affording $\theta_i$.
We have
\begin{align*}
\dim V_3(\rho_1) = 9, \qquad \dim V_3(\rho_2) = \dim V_3(\rho_1, \chi_2) = 5, \\
\dim V_3(\rho_3) = 6, \qquad \dim V_3(\rho_4) = \dim V_3(\rho_3, \chi_2) = 2,
\end{align*} 
and $\rho_4$ is not suitable. Taking projective coordinates 
$[x_1 : x_2 : x_3 : x_4 : y_1 : y_2]$, the representation $\rho_1$ is given by 
\begin{align*}
(1234)[x_1 : x_2 : x_3 : x_4 : y_1 : y_2] = [x_2 : x_3 : x_4 : x_1 :y_1 :-y_2], \\
(12)[x_1 : x_2 : x_3 : x_4 : y_1 : y_2] = [x_2 : x_1 : x_3 : x_4 : y_1 : -y_2], 
\end{align*}
and it satisfies symplectic conditions. The linear space $V_3(\rho_1)$ is 
a $\rho_1((12))$-invariant subspace of (\ref{eqA4-1}):
\begin{multline} \label{eqS4-1}
\mathrm{Span} \{ s_1(x)^3, \ s_1(x)^2y_1, \ s_1(x)y_1^2, \ s_1(x)y_2^2, \ y_1^3, \ y_1y_2^2, \\ 
s_1(x)s_2(x), \ s_2(x) y_1, \ s_3(x) \}.
\end{multline}
A general member of $V_3(\rho_1)$ is smooth since it contains a smooth cubic 
\[
s_3(x) + y_1^2y_2 + y_2^3 = 0.
\] 
We have $\dim \C_{\rho_1}(\Sy_4) = 2^2 + 1^2 + 1^2 = 6$ and 
\[
\dim V_3(\rho_1) - \dim \C_{\rho_1}(\Sy_4) = 3 = 20 - r(\Sy_4).
\]
Since we have 
\[
\C_{\rho_2}(\Sy_4) = \C_{\rho_1 \otimes \chi_1}(\Sy_4) = \C_{\rho_1}(\Sy_4) 
\]
and $\dim V_3(\rho_2) = 5$, the linear space $V_3(\rho_2)$ does not contain a smooth cubic.
\\ \indent
The representation $\rho_3$ is given by 
\begin{align*}
(1234)[x_1 : x_2 : x_3 : x_4 : y_1 : y_2]  = [x_2 : x_3 : x_4 : x_1 : -\w y_2 : -\w^2 y_1],  \\
(12)[x_1 : x_2 : x_3 : x_4 : y_1 : y_2]  = [x_2 : x_1 : x_3 : x_4 : -\w^2 y_2 : -\w y_1]. 
\end{align*}
and it satisfies the symplectic conditions. Moreover the restriction of this action to 
$\A_4$ coincides with the action (\ref{ACT-3.A4}).
The linear space $V_3(\rho_3)$ is a $\rho_3((12))$-invariant subspace of (\ref{eqA4-2})
\begin{multline} \label{eqS4-2}
{\rm Span} \{ s_1(x)^3, \quad s_1(x)s_2(x), \quad s_3(x), \quad s_1(x)y_1y_2, \\
y_1^3 + y_2^3, \quad \w^2 f(x)y_2 + \w g(x)y_1 \}.
\end{multline}
A general member of $V_3(\rho_3)$ is smooth since it contains a smooth cubic 
\[
s_3(x) + y_1^2y_2 + y_2^3 = 0.
\] 
We have $\dim \C_{\rho_3}(\Sy_4) = 1^2 + 1^2 + 1^2 = 3$ and 
\[
\dim V_3(\rho_3) - \dim \C_{\rho_3}(\Sy_4) = 3 = 20 - r(\Sy_4).
\]

%% file: r18.tex
\section{Groups of $r(G) = 18$}
We have eight groups:
\[
3^{1+4}_+:2, \quad \A_{4,3}, \quad \A_5, \quad 3^2.4, \quad \Sy_{3,3}, \quad 
\mathrm{F}_{21}, \quad \mathrm{Hol}(5), \quad \mathrm{QD}_{16}. 
\]
\subsection{The group $\mathbf{3^{1+4}_+:2, \ F_{21}}$ and  $\mathbf{QD_{16}}$} 
\label{SS-r18}
These groups are studied in \cite{LZ22}.
\\
(i) \ $G \cong 3^{1+4}_+:2$. For an appropriate choice of coordinates, the defining equations
of the corresponding cubic fourfolds are given by
\begin{align} \label{eq3^{1+4}:2}
F_1(x_1, x_2, x_3) + F_2(x_4, x_5, x_6) = 0
\end{align}
An element of order 3 in $G$ is $\diag(1, 1, 1, \w, \w, \w)$. Moreover, any smooth cubic 
fourfold with a symplectic automorphism which can be diagonalized as $\diag(1, 1, 1, \w, \w, \w)$
 has this form, and a generic such cubic fourfold has symplectic automorphism group $3^{1+4}_+:2$.
\begin{rem} \label{extraspecial}
Extraspecial groups $G = p^{1+2n}_{\pm}$ are non-abelian $p$-groups with center 
$Z(G) \cong \C_p$ and $G/Z(G) \cong \C_p^{2n}$. The exponent of $p^{1+2n}_+$ 
is $1$, and that of $p^{1+2n}_-$ is $2$. We give an explicit matrix representation 
of $3^{1+4}_+2$ in section \ref{SS-3^1+4:2.2}. 
\end{rem}
\noindent
(ii) \ $G \cong \mathrm{F}_{21}$. For an appropriate choice of coordinates, 
the defining equations of the corresponding cubic fourfolds are given by
\begin{align} \label{eqF21}
x_1^2 x_2 + x_2^2 x_3 + x_3^2 x_4 + x_4^2 x_5 + x_5^2 x_6 + x_6^2 x_1 + 
a x_1 x_3 x_5 + b x_2 x_4 x_6 = 0.
\end{align}
The automorphisms $g_1 = \diag(\z_7, \z_7^5, \z_7^4, \z_7^6, \z_7^2, \z_7^3)$ and 
$g_2 : x_i \rightarrow x_{i+2}$ generate $\mathrm{F}_{21} \cong 7.3$.
Moreover, any smooth cubic fourfold with symplectic automorphism of order 7 has
this form, and a generic such cubic fourfold has symplectic automorphism group $\mathrm{F}_{21}$.
\\ \\
(iii) \ $G \cong \mathrm{QD}_{16}$. For an appropriate choice of coordinates, the defining 
equations of the corresponding cubic fourfolds belong to
\begin{align} \label{eqQD16}
\mathrm{Span}\{ x_1^3, \ x_1x_2^2, \ x_2x_3^2, \ x_2x_4^2, \ x_1x_3x_4,\ x_4x_5^2,\ 
x_3x_6^2, \ x_2x_5x_6 \}
\end{align}
An element of order 8 in $G$ is $a = \diag(1, -1, i, -i, \z_8, \z_8^3)$. 
Moreover, any smooth cubic fourfold with a symplectic automorphism of order 8 has 
this form, and a generic such cubic fourfold has symplectic automorphism group 
$\mathrm{QD}_{16}$. 
\begin{rem} 
The quasi-dihedral group $\mathrm{QD_{16}}$ is defined by
\[
\mathrm{QD_{16}} = \left< a, b \ | \ a^8 = b^2 = 1, \ bab = a^3 \right>.
\]
A general member in (\ref{eqQD16}) is normalized into
\begin{multline}
a_1 x_1^3 + a_2 x_1^2 x_3 + a_3 (x_2 x_3^2 + x_2 x_4^2) + a_3 x_1 x_3 x_4 \\
+ a_4(x_4 x_5^2 + x_3 x_6^2) + a_5 x_2 x_5 x_6
\end{multline}
by diagonal actions of $(\CC^{\times})^6$. Then $\mathrm{QD}_{16}$ is generated by 
$a = \diag(1, -1, i, -i, \z_8, \z_8^3)$ and a permutation of coordinates by $(34)(56)$.
\end{rem}
\subsection{The group $\A_{4,3} = (\Sy_4 \times \Sy_3) \cap \A_7$}
The Schur multiplier of $\A_{4,3}$ is $\C_6$, and we must consider a covering group 
$3.\A_{4,3}$. Firstly, we consider linear actions of $\A_{4,3}$ itself. 
Tthe character table of $\A_{4,3}$ is as follows:
\[
 \begin{array}{c|rrrrrrrrr}
  \rm class&\rm1&\rm2A&\rm2B&\rm3A&\rm3B&\rm3C&\rm3D&\rm4A&\rm6A\cr
  \rm size&1&3&18&2&8&8&8&18&6\cr
\hline
  \chi_{1}&1&1&1&1&1&1&1&1&1\cr
  \chi_{2}&1&1&-1&1&1&1&1&-1&1\cr
  \chi_{3}&2&2&0&2&-1&-1&-1&0&2\cr
  \chi_{4}&2&2&0&-1&2&-1&-1&0&-1\cr
  \chi_{5}&2&2&0&-1&-1&-1&2&0&-1\cr
  \chi_{6}&2&2&0&-1&-1&2&-1&0&-1\cr
  \chi_{7}&3&-1&-1&3&0&0&0&1&-1\cr
  \chi_{8}&3&-1&1&3&0&0&0&-1&-1\cr
  \chi_{9}&6&-2&0&-3&0&0&0&0&1\cr
\end{array}
\quad 
\begin{array}{l}
\rm 2A : (12)(34) \\
\rm 2B : (12(56) \\
\rm 3A : (567) \\
\rm 3B : (123) \\
\rm 3C : (123)(567) \\
\rm 3D : (123)(576) \\
\rm 4A : (1234)(67) \\
\rm 6A : (12)(34)(567)
\end{array}
\]
Since we have $\widehat{\A_{4,3}} = \C_2$, we do not have to consider semi-invariants. 
A character $\phi$ satisfying conditions 
\[
|\phi(2A)| = |\phi(2B)| =2, \qquad |\phi(4A)| =0
\] 
is one of the followings:
\[
\phi_i = \chi_1 + \chi_i + \chi_8, \quad 
\phi_i \otimes \chi_2 = \chi_2 + \chi_i + \chi_7 \quad (i = 4,5,6). 
\]
We have $\Out(\A_{4,3}) \cong \Sy_3$ which permutes classes 3B, 3C and 3D, and 
hence $\phi_4, \ \phi_5$ and $\phi_6$.  Hence we consider only $\phi_4$ and 
 $\phi_4 \otimes \chi_2$. 
\begin{rem}
The conjugation by $(67) \in \Sy_7$ interchages classes {\rm 3C} and {\rm 3D}.  
However, the cycle type of {\rm 3B} is different with that of {\rm 3B} and {\rm 3C}, 
and we need another kind of automorphisms to generate $\Out(\A_{4,3})$. 
One of them is given by
\begin{multline*}
[ (3,4)(6,7),\ (5,6,7), \ (2,3,4), \ (1,4)(2,3), \ (1,3)(2,4) ] \\
 \longrightarrow
 [(2,4)(5,6), \ (5,6,7), \ (2,4,3)(5,6,7), \ (1,4)(2,3), \ (1,2)(3,4)],
\end{multline*}
which is obtained by GAP.
\end{rem}
Let $\rho$ be a representation affording $\phi_4$. 
Since Molien series of $\rho$ and $\rho \otimes \chi_2$ are
\begin{align*}
\mu(\rho, t) = 1 + t + 3t^2 + 5t^3 + \cdots, \qquad
\mu(\rho \otimes \chi_2, t) = 1 + 3t^2 + t^3 + \cdots,
\end{align*}
we have $\dim V_3(\rho \otimes \chi_2) = 1$ and this is not the case. 
The representaion $\rho$ is given as follows. 
Let $[x_1 : x_2 : x_3 : x_4 : y_1 : y_2 : y_3]$ be projective coordinates of $\PP^6$, 
and $H_0$ be a hyperplane
\[
s_1(x) + s_1(y) = (x_1 + x_2 + x_3 + x_4) + (y_1 + y_2 + y_3) = 0
\]
Now $\A_{4,3}$ acts on $H_0 \cong \PP^5$ as permutations of 
$x_i$ and $y_i$ respectively, and we have the following decomposition
\begin{align*}
&\chi_1 \text{-space} \ [z : z : z : z : -z : -z : -z], \\ 
&\chi_8 \text{-space} \ [x_1 : x_2 : x_3 : x_4 : 0 : 0 : 0] \ \text{with} \ s_1(x)=0,  \\
&\chi_4 \text{-space} \ [0 : 0 : 0 : 0 : y_1 : y_2 : y_3] \ \text{with} \ s_1(y)=0.
\end{align*}
The linear space $V_3(\rho)$ in the quotient ring $\CC[x_i, y_i]/(s_1(x)+s_1(y))$ is 
given by 
\begin{align} \label{eqA43-1}
\mathrm{Span}\{ s_1(x)^3, \ s_1(x)s_2(x), \ s_1(x)s_2(y), \ s_3(x), \ s_3(y) \}.
\end{align}
Note that this family contains the Clebsch-Segre cubic (\ref{Clebsch-Segre}), and a 
general member is smooth. 
We have $\dim \C_{\rho}(\A_{4,3}) = 1^2 + 1^2 + 1^2 = 3$ and 
\[
\dim V_3(\rho) - \dim \C_{\rho}(\A_{4,3}) = 2 = 20 - r(\A_{4,3}).
\]
\indent
Next we consider a triple cover $3.A_{4,3} = (\A_{4,3})^*/\C_2$ whose GAP-ID is (216,95). 
Since we have $\widehat{G^*} \cong \widehat{G} \cong \C_2$, 
 we do not have to consider semi-invariants. Irreducible representaions of $G^*$ 
with a non-trivial action of the Schur multiplier are as follows.
\begin{align*}
\begin{array}{l|rrrrrrrrrrr}
  \rm class&\rm1 & \rm2A & \rm2B &
  \rm3A & \rm3B & \rm3C & \rm3D & \rm3E & \rm3F & \rm4 & \cdots \cr
  \rm size&1&3&18&1&1&6&24&24&24&18 & \cdots \cr
\hline
  \chi_{10} &3&3&1&3\w&3\bw&0&0&0&0&1 \cr
  \chi_{11} & 3&-1&1&3\w&3\bw&0&0&0&0&-1 \cr
  \chi_{12} = \chi_{10} \otimes \chi_2 &3&3&-1&3\w&3\bw&0&0&0&0&-1 \cr
  \chi_{13} = \chi_{11} \otimes \chi_2 &3&-1&-1&3\w&3\bw&0&0&0&0&1 \cr
  \chi_{14} = \overline{\chi_{10}} &3&3&1&3\bw&3\w&0&0&0&0&1 \cr
  \chi_{15} = \overline{\chi_{11}} &3&-1&1&3\bw&3\w&0&0&0&0&-1 \cr 
  \chi_{16} = \overline{\chi_{10}} \otimes \chi_2 &3&3&-1&3\bw&3\w&0&0&0&0&-1 \cr
  \chi_{17} = \overline{\chi_{11}} \otimes \chi_2 &3&-1&-1&3\bw&3\w&0&0&0&0&1 \cr
  \chi_{18} &6&-2&0&6\w&6\bw&0&0&0&0&0 \cr
  \chi_{19} &6&-2&0&6\bw&6\w&0&0&0&0&0 \cr
\end{array}
\\
\text{(Conjugacy classes for order $6$ and $12$ are omitted.)}
\end{align*}
Classes 1, 3A and 3B form the Schur multiplier, and they act as scalar matrices for
\[
\chi_i + \chi_j \ (10 \leq i,j \leq13), \quad
\chi_k + \chi_l \ (14 \leq i,j \leq17), \quad
\chi_{18}, \quad \chi_{19}. 
\]
Among these characters, conditions
\[
|\phi(2A)| = |\phi(2B)| = 2, \qquad |\phi(4A)|=0
\]
are satisfied only by
\[
\phi_1 = \chi_{10} + \chi_{11}, \quad \overline{\phi_1}, \quad
\phi_2 = \chi_{12} + \chi_{13}, \quad \overline{\phi_2}.
\]
Let $\rho_i$ be a representation affording $\phi_i$. We have
Molien series
\[
\mu(\rho_1, t) = 1 + 4t^3 + \cdots, \qquad
\mu(\rho_2, t) = 1 + 17t^6 + \cdots
\]
and
\[
\dim V_3(\rho_1) = 4, \qquad \dim V_3(\rho_2) = 0.
\]
For $a$ and $b$ in (\ref{ab in 3.A7}), the following elements
\begin{align}
u_1 = aba^{-1}bab^{-2}a^{-1}b, \quad
u_2 = b^2ab^2a^{-1}bab, \quad
u_3 = b^2aba^{-1}b^{-1}ab^2
\end{align}
generate a group with GAP-ID (216, 95), and
this gives one of $\rho_1$ and $\overline{\rho_1}$. Moreover we have
\[
 \overline{u_1} = u_1, \quad \overline{u_2} = u_3^{-1}(u_2 u_1)^2, \quad
\overline{u_3} = u_3u_2^2u_3u_2,
\]
and hence $V_3(\rho_1) = V_3(\overline{\rho_1})$. This space is given by
\begin{align} \label{eqA43-2}
\mathrm{Span} \{ p_1(x),\ p_2(x),\ p_3(x),\ p_4(x)\}
\end{align}
where
\begin{align*}
p_1(x) &= x_2^3 + x_4^3 - 3 x_{246} + x_6^3 \\
p_2(x) &= 8 x_{135} + 4 x_{235} + 4 x_{145} + 2 x_{245} + 4 x_{136} + 
 2 x_{236} + 2 x_{146} + x_{246} \\
p_3(x) &= 2 x_{122} + 2 x_{344} + 4 x_{245} + 4 x_{236} + 4 x_{146} + 9 x_{246} + 2 x_{566} \\
p_4(x) &= 2 x_1^3 + 3 x_{112} + 2 x_3^3 + 3 x_{334} - 3 x_{245} + 2 x_5^3 - 
 3 x_{236} - 3 x_{146} - 6 x_{246} + 3 x_{556}
\end{align*}
with $x_{ijk} = x_i x_j x_k$. This family contains a smooth cubic fourfold 
(\ref{3.A7-invariant}), and a general member is smooth. 
We have $\dim \C_{\rho_1}(3.\A_{4,3}) = 1^2 + 1^2 = 2$ and 
\[
\dim V_3(\rho_1) - \dim \C_{\rho_1}(3.\A_{4,3}) = 2 = 20 - r(\A_5).
\]
\subsection{The alternating group $\mathbf{\A_5}$} \label{SS-A5}
The Schur multiplier of $\A_5$ is $\C_2$, and we do not have to consider 
a Schur cover. The character tables of $\A_5$ is as follows:
\[
\begin{array}{c|rrrrr}
  \rm class&\rm1&(12)(34) & (123) & (12345) &(13524) \cr
  \rm size&1&15 \ &20 \ &12 \ &12 \ \cr
\hline
  \chi_1 &1&1 \ &1 \ &1 \ &1 \ \cr
  \chi_2 &3&-1 \ &0 \ &\frac{1+\sqrt{5}}{2} \ &\frac{1-\sqrt{5}}{2} \ \cr
  \chi_3 &3&-1 \ &0 \ &\frac{1-\sqrt{5}}{2} \ &\frac{1+\sqrt{5}}{2} \ \cr
  \chi_4 &4&0 \ &1 \ &-1 \ &-1 \ \cr
  \chi_5 &5&1 \ &-1 \ &0 \ &0 \ \cr
\end{array}
\]
Since we have $\widehat{\A_5} = 1$, we do not have to consider semi-invariants. 
If $\rho : \A_5 \rightarrow \GL_6(\CC)$ satisfies the symplectic condition, we have
\[
|\Tr_{\rho}(12)(34)| =2, \quad |\Tr_{\rho}(123)| = 0,3, \quad 
|\Tr_{\rho}(12345)| = 1.
\]
These conditions are satisfied only for
\begin{align} \label{theta-A5}
\theta_1 = 2 \chi_1 + \chi_4, \quad \theta_2 = \chi_1 + \chi_5, \quad
\theta_3 = \chi_2 + \chi_3.
\end{align}
For representations $\rho_i$ affording $\theta_i$, we have Molien series
\begin{align*}
\mu(\rho_1, t) &=  1 + 2t + 4t^2 + 7t^3 + \cdots, \\
\mu(\rho_2, t) &=  1 + t + 2t^2 + 4t^3 + \cdots, \\
\mu(\rho_3, t) &=  1 + 2t^2 + 6t^4 + \cdots, 
\end{align*}
and $V_3(\rho_3) = 0$. For projective coordinates $[x_1: \cdots : x_5 : y]$, 
the representation $\rho_1$  is given as permutations of $x_i$, and we have 
\begin{multline} \label{eqA5-1}
V_3(\rho_1) = \mathrm{Span} \{
 y^3, \ y^2 s_1(x), \ y s_1(x)^2, \ y s_2(x), \ s_1(s)^3, \ 
 s_1(x) s_2(x), \ s_3(x) \}
\end{multline}
Although $F \in V_3(\rho_1)$ is stable under permutatoins by $\Sy_5$, the action of 
$(12) \in \Sy_5$ is non-symplectic. 
Replacing $x_5$ and $y$ by $y_1$ and $y_2$, we see that this is a subspace of 
(\ref{eqA4-1}).  
We have $\dim \C_{\rho_1}(\A_5) = 2^2 + 1^2 = 5$ and 
\[
\dim V_3(\rho_1) - \dim \C_{\rho_1}(\A_5) = 2 = 20 - r(\A_5).
\]
The representation $\rho_2$ is given by permutations of projective coordinates 
$[x_1: \cdots : x_6]$ by
``non-standard'' $\A_5 \subset \A_6$ generated by $(12345)$ and $(35)(46)$
(see Appendix A). Putting $x_{ijk} = x_i x_j x_k$ and
\begin{align*}
I_1 &= x_{124} + x_{134} + x_{135} + x_{235} + x_{245} + x_{126} + 
 x_{236} + x_{346} + x_{156} + x_{456}, \\
I_2 &= x_{123} + x_{234} + x_{125} + x_{145} + x_{345} + x_{136} + 
 x_{146} + x_{246} + x_{256} + x_{356},
\end{align*}
we have
\begin{align} \label{eqA5-2}
V_3(\rho_2) = \mathrm{Span} \{
s_1(x)^3, \quad s_1(x)s_2(x), \quad s_3(x), \quad I_1, \quad I_2 \}
\end{align}
with a linear relaton 
\[
s_1^3 - 6 (I_1 + I_2) - 3 s_1 s_2 + 2 s_3 = 0. 
\]
A general $F \in V_3(\rho_2)$ defines a smooth cubic fourfold since $s_3(x) \in V_3(\rho_2)$ 
defines the Fermat cubic. 
We have $\C_{\rho_2}(\A_5) = 1^2 + 1^2 = 2$ and
\[
\dim V_3(\rho_2) - \dim \C_{\rho_2}(\A_5) = 2 = 20 - r(\A_5)
\]
as desired.
\subsection{The group $\mathbf{G \cong 3^2.4}$ with GAP-ID (36,9)} 
\label{SS-3^2.4}
This group $G$ is represented as a (maximal) subgroup of $\A_6$ generated by
$(123), \ (456), \ (1425)(36)$, and containing $\A_{3,3} = 3^2:2$ as a subgroup of index 2. 
We have a Schur cover $G^* \cong 3.G$ with GAP-ID (108,15).
Firstly, we consider linear representations of $G$ itself. The character table of $G$ is 
as follows.
\[
\begin{array}{c|rrrrrr}
  \rm class& 1 &\rm 2A &\rm 3A & \rm 3B &\rm 4A &\rm 4B \cr
  \rm size&1&9&4&4&9&9\cr
\hline
  \eta_{1}&1&1&1&1&1&1\cr
  \eta_{2}&1&1&1&1&-1&-1\cr
  \eta_{3}&1&-1&1&1&i&-i\cr
  \eta_{4}&1&-1&1&1&-i&i\cr
  \eta_{5}&4&0&1&-2&0&0\cr
  \eta_{6}&4&0&-2&1&0&0\cr
\end{array}
\qquad
\begin{array}{l}
\rm 2A : (12)(45) \\
\rm 3A : (123) \\
\rm 3B : (123)(456) \\
\rm 4A : (1425)(36) \\
\rm 4B : (1524)(36)
\end{array}
\]
We have $\Out(G) \cong \C_2 \times \C_2$, and it interchanges classes 
nA with nB for $n = 3,4$.  
Since we have $\widehat{G} = \C_4$, we consider only invariants with the trivial character. 
Faithful characters satisfying the condition on absolute values are
\begin{align*}
 \phi_1 = \eta_1 + \eta_2 + \eta_5, \quad \phi_2 = \eta_1 + \eta_2 + \eta_6, \\
 \phi_1 \otimes \eta_3 = \eta_3 + \eta_4 + \eta_5, \quad 
 \phi_2 \otimes \eta_3 =\eta_3 + \eta_4 + \eta_6.
\end{align*}
Since $\Out(G)$ interchanges $\eta_5$ with $\eta_6$, we consider only 
$\phi_1$ and $\phi_1 \otimes \eta_3$. 
Note that
\begin{align*}
\phi_1 | \A_{3,3} &= (\eta_1 + \eta_2 + \eta_5)| \A_{3,3} 
= (2\chi_1 + \chi_3 + \chi_4) = \theta_{34}, \\
\phi_1 \otimes \eta_3 | \A_{3,3} &= \theta_{34} \otimes \chi_2
\end{align*}
for $\chi_i$ and $\theta_{34}$ in section \ref{SS-A33}. A representaion $\rho$ 
affording $\phi_1$ is given by permutatioins of coordinates $[x_1 : x_2 : x_3 : y_1 : y_2 : y_3]$, 
and $V_3(\rho)$ is $(1425)(36)$-invariants of (\ref{eqA33-1}):
\begin{multline} \label{eq3^2.4-1}
\mathrm{Span} \{
s_1(x)^3 + s_1(y)^3, \quad s_2(x)s_1(x) + s_2(y)s_1(y), \quad s_3(x) + s_3(y), \\
s_1(x)s_2(y) + s_1(y)s_2(x), \quad s_1(x)s_1(y)^2 + s_1(x)^2 s_1(y) \}.
\end{multline}
We have $\dim \C_{\rho}(G) = 1^2 + 1^2 + 1^2 = 3$ and 
\[
\dim V_3(\rho) - \dim \C_{\rho}(G) = 2 = 20 - r(G).
\]
For a representaion $\rho \otimes \eta_3$, there is no smooth cubic since  
$V(\rho \otimes \eta_3)$ is a subspace of (\ref{eqA33-singular}). 
\\ \indent
Next, we consider $G^* = 3.G$. Since we have $\widehat{G^*} \cong \widehat{G} \cong \C_4$, 
 we do not have to consider semi-invariants. Irreducible representaions of $G^*$ 
with a non-trivial action of the Schur multiplier are given as follows.
\[
\begin{array}{c|rrrrrrrrrrrrrr}
  \rm class&\rm1&\rm2&\rm3A&\rm3B&\rm3C&\rm3D&\rm4A&\rm4B&\rm6A&\rm6B&
  \rm12A&\rm12B&\rm12C&\rm12D\cr
  \rm size&1&9&1&1&12&12&9&9&9&9&9&9&9&9\cr
\hline
  \eta_{7}&3&1&\al&\ba&0&0&i&-i&\w&\bw&-\de&\bd&-\bd&\de\cr
  \eta_{8}&3&1&\al&\ba&0&0&-i&i&\w&\bw&\de&-\bd&\bd&-\de\cr
  \eta_{9}&3&-1&\al&\ba&0&0&1&1&-\w&-\bw&\w&\bw&\bw&\w\cr
  \eta_{10}&3&-1&\al&\ba&0&0&-1&-1&-\w&-\bw&-\w&-\bw&-\bw&-\w\cr 
  \eta_{11}&3&1&\ba&\al&0&0&-i&i&\bw&\w&-\bd&\de&-\de&\bd\cr
  \eta_{12}&3&1&\ba&\al&0&0&i&-i&\bw&\w&\bd&-\de&\de&-\bd\cr
  \eta_{13}&3&-1&\ba&\al&0&0&1&1&-\bw&-\w&\bw&\w&\w&\bw\cr
  \eta_{14}&3&-1&\ba&\al&0&0&-1&-1&-\bw&-\w&-\bw&-\w&-\w&-\bw\cr
\end{array}
\]
\[
\al = 3 \w, \quad \de = i \w, \quad \bd = -i\bw
\]
Characters $\eta_{11}, \cdots, \eta_{14}$ are complex conjugates of $\eta_7, \cdots, \eta_{10}$.
Classes 1, 3A and 3B form the Schur multiplier $\C_3$ which act as scalar matrices for
\[
\eta_i + \eta_j \quad (7 \leq i,j \leq 10), \qquad
\eta_k + \eta_l \quad (11 \leq k,l \leq 14).
\]
The condition $(\eta_i + \eta_j)(4A) = 0$ is satisfied only by
\[
\eta_7 + \eta_8, \quad \eta_9 + \eta_{10}, \quad \eta_{11} + \eta_{12}, \quad 
\eta_{13} + \eta_{14}.
\]
For representations $\rho_i$ affording $\eta_i$, we have Molien series
\begin{align*}
 \mu(\rho_i + \rho_j, \, t) = \sum_{d=0}^{\infty} (\dim V_d(\rho_i + \rho_j))t^d 
= 1 + 4 t^3 + \cdots \qquad ((i,j) = (7,8),(11,12)) \\
 \mu(\rho_i + \rho_j, \, t) = \sum_{d=0}^{\infty} (\dim V_d(\rho_i + \rho_j))t^d 
= 1 + 20 t^6 + \cdots \qquad ((i,j) = (9,10),(13,14)), 
\end{align*}
and there is no invariant cubic for $\rho_9 + \rho_{10}$ and $\rho_{13} + \rho_{14}$. 
\\ \indent
Let us consider a matrix $\tilde{h} = h \oplus (-h)$ where
\[
h = \frac{1}{\sqrt{3}} \begin{bmatrix} 
1 & 1 & 1 \\ 1 & \w & \w^2 \\ 1 & \w^2 & \w \end{bmatrix}.
\] 
Eigenvalues of $h$ are $1,-1,i$, and hence $\tilde{h}$ satisfies the symplectic 
condition. 
The group $H \cong 3.A_{3,3}$ generated by $\tilde{a}, \tilde{b}, \tilde{c}$ 
in (\ref{abc2-A33}) is normalized by $\tilde{h}$, and $H$ together with 
$\tilde{h}$ generate a group $H_2 \cong 3.(3^2.4)$ with GAP-ID (108,15) 
(we have $\tilde{h}^2 = \tilde{a}$ and $[H_2:H]=2$).  This gives one of 
$\rho_7 + \rho_8$ and $\rho_{11} + \rho_{13}$. Moreover we have $\overline{H} = H$ 
(see (\ref{abc-conj})) and $\overline{h} = h^3$, that is, $\overline{H_2} = H_2$. 
Therefore we have
\[
V_3(\rho_7 + \rho_8) = V_3(\overline{\rho_7} + \overline{\rho_8}) 
= V_3(\rho_{11} + \rho_{13})
\] 
The space $V_3(\rho_7 + \rho_8)$ is $\tilde{h}$-invariants of
(\ref{eqA33-2}):
\begin{align} \label{eq3^2.4-2}
\mathrm{Span} \{ s_3(x) - 3u_- x_1 x_2 x_3, \quad s_3(y) - 3u_+ y_1 y_2 y_3, \quad
I_1 - u_- I_2, \quad I_3 - u_+ I_4 \}
\end{align}
where $u_{\pm} = 1 \pm \sqrt{3}$. 
We have $\dim \C_{\rho_7 + \rho_8}(G^*) = 1^2 + 1^2 = 2$ and 
\[
\dim V_3(\rho_7 + \rho_8) - \dim \C_{\rho_7 + \rho_8}(G^*) = 2 = 20 - r(G).
\]
\subsection{The groups $\mathbf{\Sy_{3,3} = \Sy_3 \times \Sy_3}$} 
\label{SS-S33}
The Schur multiplier of $\Sy_{3,3}$ is $\C_2$, and we do not have to consider a Schur cover. 
The character table of $\Sy_3$ is 
\[
\begin{array}{c|rrr}
  \rm class&\rm1&(12)&(123)\cr
  \rm size&1&3&2\cr
\hline
  \eta_{1}&1&1&1\cr
  \eta_{2}&1&-1&1\cr
  \eta_{3}&2&0&-1\cr
\end{array}
\]
and irreducible characters of $\Sy_{3,3}$ are given as tensor products 
$\eta_{ij} = \eta_i \boxtimes \eta_j$ where
\[
\eta_{ij}(a, b) = \eta_i (a)\eta_j(b) \qquad (a,b) \in \Sy_3 \times \Sy_3,
\] 
and the character table of $\Sy_{3,3}$ is as follows.
\begin{align*}
\begin{array}{c|rrrrrrrrr}
  \rm class&\rm1&\rm2A&\rm2B&\rm2C&\rm3A&\rm3B&\rm3C&\rm6A&\rm6B\cr
  \rm size&1&3&3&9&2&2&4&6&6\cr
\hline
  \eta_{11}&1&1&1&1&1&1&1&1&1\cr
  \eta_{21}&1&-1&1&-1&1&1&1&1&-1\cr
  \eta_{12}&1&1&-1&-1&1&1&1&-1&1\cr
  \eta_{22}&1&-1&-1&1&1&1&1&-1&-1\cr
  \eta_{13}&2&2&0&0&2&-1&-1&0&-1\cr
  \eta_{23}&2&-2&0&0&2&-1&-1&0&1\cr
  \eta_{31}&2&0&2&0&-1&2&-1&-1&0\cr
  \eta_{32}&2&0&-2&0&-1&2&-1&1&0\cr
  \eta_{33}&4&0&0&0&-2&-2&1&0&0\cr
\end{array}
\qquad
\begin{array}{l}
\rm 2A = (12) \\
\rm 2B = (45) \\
\rm 2C = (12)(45) \\
\rm 3A = (123) \\
\rm 3B = (456) \\
\rm 3C = (123)(456) \\
\rm 6A = (123)(45) \\
\rm 6B = (12)(456) \\
\end{array}
\end{align*}
Since we have $\widehat{\Sy_{3,3}} \cong C_2 \times \C_2$,we consider only 
invariants with the trivial character. We see that $\Sy_{3,3}$-invariant cubics are 
given as subspaces of $\A_{3,3}$-invariant cubics $V_3(\rho)$ in (\ref{eqA33-1}).  
We have
\begin{align*}
\eta_{11}|\A_{3,3} &= \eta_{22}|\A_{3,3} = \chi_1, \quad
\eta_{21}|\A_{3,3} = \eta_{12}|\A_{3,3} = \chi_2, \\
\eta_{13}|\A_{3,3} &= \eta_{23}|\A_{3,3} = \chi_3, \quad
\eta_{31}|\A_{3,3} = \eta_{32}|\A_{3,3} = \chi_4, \quad
\eta_{33}|\A_{3,3} = \chi_5 + \chi_6. 
\end{align*}
for $\chi_i$ in section \ref{SS-A33}, and a character $\phi$ satisfying   
$|\phi(2A)| = 2$ and 
\[
\phi|\A_{3,3} = \theta_{34} = 2\chi_1 + \chi_3 + \chi_4
\]
is one of the followings:
\begin{align*}
\phi_1 = \eta_{11} + \eta_{22} + \eta_{13} + \eta_{31}, \quad 
\phi_1 \otimes \eta_{22} = \eta_{11} + \eta_{22} + \eta_{23} + \eta_{32}, \\
\phi_2 = \eta_{11} + \eta_{22} + \eta_{13} + \eta_{32}, \quad
\phi_2 \otimes \eta_{22} = \eta_{11} + \eta_{22} + \eta_{23} + \eta_{31}.
\end{align*}
Let us consider a representation $\rho_1$ affording $\phi_1$.  
Permutations of coordinates $[x_1:x_2:x_3:y_1:y_2:y_3]$ by $\A_{3,3}$ is normalized by
an action of a matrix
\begin{align*}
b_1 = \frac{1}{3}\begin{bmatrix} -2 & 1 & -2 \\ 1 & -2 & -2 \\ -2 & -2 & 1 \end{bmatrix}
\oplus \diag(1,1,1)
\end{align*}
with eigenvalues $(-1,-1,1,1,1,1)$.  These actions generate $\Sy_{3,3}$ and this gives $\rho_1$.  
By the involution $b_1$, the space $V_3(\rho)$ in (\ref{eqA33-1}) is decomposed into 
the $(+1)$-eigenspace
\begin{multline} \label{eqS33}
V_3(\rho_1) = \mathrm{Span} \{ s_1(y)^3, \ s_2(y)s_1(y), \ s_3(y), \ s_1(x)^2 s_1(y), \ 
s_2(x) s_1(y), \ p(x) \} 
\end{multline}
where $p(x) = 2s_1(x)^3 -9 s_1(x)s_2(x) + 9s_3(x)$, and the $(-1)$-eigenspace
\begin{multline} 
 V_3(\rho_1, \eta_{22}) = V_3(\rho_1 \otimes \eta_{22})  \\
=\mathrm{Span} \{ s_1(x)^3, \ s_2(x)s_1(x), \ s_1(x)s_2(y), \ s_1(x)s_1(y)^2 \}.
\end{multline}
A general member of $V_3(\rho_1)$ defines a smooth cubic fourfold since 
\[
s_3(y) + s_2(x)s_1(y) + p(x) = 0
\]
is smooth. On the other hand, $F \in V_3(\rho_1 \otimes \eta_{22})$ is reducible. 
We have 
\[
\dim \C_{\rho_1}(\Sy_{3,3}) = 1^2 + 1^2 + 1^2 + 1^2 = 4
\]
and
\[
\dim V_3(\rho_1) - \dim \C_{\rho_1}(\Sy_{3,3}) = 2 = 20 - r(\Sy_{3,3}).
\] 
Similarly,  $A_{3,3}$ and an involution
\begin{align*}
b_2 = \begin{bmatrix} 0 & -1 & 0 \\ -1 & 0 & 0 \\ 0 & 0 & -1 \end{bmatrix}
\oplus \diag(1,1,1)
\end{align*}
generate $\Sy_{3,3}$, which gives a representation $\rho_2$ affording $\phi_2$.   
We can decompose $V_3(\rho)$ into $(+1)$-eigenspace of $b_2$
\begin{align*} 
V_3(\rho_2) = \mathrm{Span} \{s_1(y)^3, \ s_2(y)s_1(y), \ s_3(y), \ 
s_1(y)s_2(x), \ s_1(x)^2 s_1(y) \},
\end{align*}
and $(-1)$-eigenspace of $b_2$
\begin{multline*} 
V_3(\rho_2, \eta_{22}) = V_3(\rho_2 \otimes \eta_{22}) \\
= \mathrm{Span} \{ s_1(x)^3, \ s_2(x)s_1(x), \ s_3(x), \ s_1(x)s_2(y), \ s_1(x)s_1(y)^2 \}.
\end{multline*}
For $F \in V_3(\rho_2)$, a cubic fourfold $F = 0$ is singular along 
$s_1(y) = s_3(y) = 0$. Similarly, $F \in V_3(\rho_2 \otimes \eta_{22})$ defines a singular cubic. 
\subsection{The group $\mathbf{Hol}(5)$}
The group $\mathrm{Hol}(5) \cong \D_{10}.2$ is given as a subgroup of $\Sy_5$ generated by 
$(12)(34)$ and $(2,3,4,5)$.  The Schur multiplier of $\mathrm{Hol}(5)$ is 1, and hence we have 
$\mathrm{Hol}(5)^* = \mathrm{Hol}(5)$.  The character table of $\mathrm{Hol}(5)$ is as follows:
\[
\begin{array}{c|rrrrr}
  \rm class&\rm1& (12)(34) & (2345) & (2534) & (12354) \cr
  \rm size & 1 & 5 \quad & 5 \quad & 5 \quad & 4 \quad \cr
\hline
  \chi_{1}&1&1 \quad &1 \quad &1 \quad &1 \quad \cr
  \chi_{2}&1&1 \quad &-1 \quad &-1 \quad &1 \quad \cr
  \chi_{3}&1&-1 \quad &-i \quad &i \quad &1 \quad \cr
  \chi_{4}&1&-1 \quad &i \quad &-i \quad &1 \quad \cr
  \chi_{5}&4&0 \quad &0 \quad &0 \quad &-1 \quad \cr
\end{array}
\]
Since $\widehat{\mathrm{Hol}(5)} \cong \C_4$, we consider only invariants with the 
trivial character. Faithful characters satisfying the condition on absolute values are
\[
\theta = \chi_1 + \chi_2 + \chi_5, \quad 
\theta \otimes \chi_3 = \chi_3 + \chi_4 + \chi_5.
\]
Let $\rho$ be a representation $\rho$ affording $\theta$. We have Molien series
\begin{align*}
\mu(\rho, t) = 1 + t + 3t^2 + 5t^3 + \cdots, \qquad
\mu(\rho \otimes \chi, t) = 1 + 2t^2 + t^3 + \cdots
\end{align*}
and 
\[
\dim V_3(\rho) = 5, \qquad \dim V_3(\rho \otimes \chi_3) = 1.
\]
Therefore $\rho \otimes \chi_3$ is not suitable. The representation $\rho$ is given by
\begin{align*}
(12)(34)[x_1: \cdots : x_5 : y] = [x_2 : x_1 : x_4 : x_3 : x_5 : y],
\\
(2345)[x_1: \cdots : x_5 : y] = [x_1 : x_3 : x_4 : x_5 : x_2 : -y],
\end{align*}
which satisfies symplectic conditions. Invariant cubics are  
\begin{align} \label{eqHol(5)}
\mathrm{Span} \{ s_1(x)^3, \quad s_1(x)s_2(x), \quad s_1(x)y^2, \quad s_3(x), \quad Q(x)y \}
\end{align}
where
\begin{multline}
Q(x) = -x_1 x_2 + x_1 x_3 - x_2 x_3 - x_1 x_4 + x_2 x_4 \\
+ x_3 x_4 + x_1 x_5 + x_2 x_5 -  x_3 x_5 - x_4 x_5.
\end{multline}
We have $\dim \C_{\rho}(\mathrm{Hol}(5)) = 1^2 + 1^2 + 1^2$ and
\[
\dim V_3(\rho) - \dim \C_{\rho}(\mathrm{Hol}(5)) = 2 = 20 - r(\mathrm{Hol}(5)).
\] 

%% file: r19.tex
\section{Groups of $r(G) = 19$}
We have the following seven groups
\begin{align*}
3^{1+4}:2.2, \quad \A_6, \quad \mathrm{L}_2(7), \quad \Sy_5, \quad \M_9, 
\quad \mathrm{N}_{72}, \quad \mathrm{T}_{48},
\end{align*}
and $\A_6$ was already considered in section \ref{SS-A6}.
\subsection{The group $\mathbf{3^{1+4}_+:2.2}$ with GAP-ID (972, 776)}
\label{SS-3^1+4:2.2}
This group $G$ contains $H \cong 3^{1+4}_+:2$ as a normal subgroup of index 2, 
and therefore cubic fourfolds such that $G \subset \Aut^s(X)$ are given as 
subfamilies of (\ref{eq3^{1+4}:2}) which are normalized into the Hesse normal forms 
\[
F_{\al}(x_1,x_2,x_3) = F_{\be}(y_1,y_2,y_3) \quad \text{where} \quad 
F_{\lambda}(x,y,z) = x^3 + y^3 + z^3 -27 \lambda xyz.
\]
For this normal form, $H \subset \PGL_6(\CC)$ is generated by natural actions of 
$\A_{3,3}$ and
\[
a = \diag(1,1,1,\w,\w,\w), \quad b_1 = \diag(1,1,1,1,\w,\w^2), \quad 
b_2 = \diag(1,\w,\w^2,1,1,1). 
\]
The center $Z(H)$ of $H$ is $\{1 ,a, a^2\}$, but $Z(G)$ is trivial for 
the group $G$ of GAP-ID (972, 776). A normalizer $g \in \PGL_6(\CC)$ of $H$ acts 
on $Z(H)$, and hence $g$ takes the form of 
$\begin{bmatrix} A & \0 \\ \0 & D \end{bmatrix}$ or 
$\begin{bmatrix} \0 & B \\ C & \0 \end{bmatrix}$. In the fomer case, $g$ commutes 
with $Z(H)$ and we see that $g$ and $H$ generate a group with non-trivial center. 
In the later cases, two cubic curves $F_{\al}(x) = 0$ and $F_{\be}(y) = 0$ are isomorphic 
if $g$ preserves the equation $F_{\al}(x) = F_{\be}(y)$. Therefore, cubic fourfolds 
such that $G \subset \Aut^s(X)$ are given by equations
\begin{align} \label{eq3^{1+4}:2.2}
F_{\al}(x_1,x_2,x_3) = F_{\al}(y_1,y_2,y_3) 
\end{align}
and $G$ is generated by $H$ and
\[
g = \begin{bmatrix} \0 & B \\ \1 & \0 \end{bmatrix} \quad \text{where} \quad
B = \begin{bmatrix} 0&1&0 \\ 1&0&0 \\ 0&0&1 \end{bmatrix}
\]
whose eigenvalues are $(1,1,-1,-1,i,-i)$. 
\subsection{The projective special linear group $\mathrm{L}_2(7)$} 
The group $\mathrm{L}_2(7) = \PSL_2(\mathbb{F}_7)$ is isomorphic to a maximal 
subgroup of $\A_7$ generated by $(1,2)(3,6)$ and $(1,2,3,4,5,6,7)$. 
The Schur multiplier of $\mathrm{L}_2(7)$ is $\C_2$, and we do not have to consider 
a Schur cover. The character table of 
$\mathrm{L}_2(7)$ is as follows:
\[
\begin{array}{c|rrrrrr}
  \rm class&\rm1&\rm2A&\rm3A&\rm4A&\rm7A&\rm7B\cr
  \rm size&1&21&56&42&24&24\cr
\hline
  \chi_{1}&1&1&1&1&1&1\cr
  \chi_{2}&3&-1&0&1& \frac{-1+\sqrt{-7}}{2} & \frac{-1-\sqrt{-7}}{2} \cr
  \chi_{3}&3&-1&0&1& \frac{-1-\sqrt{-7}}{2} & \frac{-1+\sqrt{-7}}{2} \cr
  \chi_{4}&6&2&0&0&-1&-1\cr
  \chi_{5}&7&-1&1&-1&0&0\cr
  \chi_{6}&8&0&-1&0&1&1\cr
\end{array}
\qquad
\begin{array}{l}
\rm 2A : (12)(36) \\
\rm 3A : (235)(476)\\
\rm 4A : (2347)(56) \\
\rm 7A : (1234567) \\
\rm 7B : (1237645)
\end{array}
\]
We see that $\chi_4$ is the unique character of degree 6 such that 
$\chi_4(4\A)=0$. The corresponding cubic fourfolds are given as 
\begin{align} \label{eqL2(7)}
\begin{cases}
x_1 + x_2 + x_3 + x_4 + x_5 + x_6 + x_7 = 0 \\
c_1( x_1^3 + x_2^3 + x_3^3 + x_4^3 + x_5^3 + x_6^3 + x_7^3) \\
\qquad \qquad + c_2(x_{124} + x_{235} +x_{346} + x_{156} + x_{137} + x_{457} + x_{267})= 0
\end{cases}
\end{align}
in $\PP^6$ where $x_{ijk} = x_i x_j x_k$. 
This family contains the Clebsch-Segre cubic (\ref{Clebsch-Segre}), 
and a general member is smooth. 
\subsection{The symmetric group $\Sy_5$}
The Schur multiplier of $\Sy_5$ is $\C_2$, and we do not have to consider a Schur cover. 
The character table of $\Sy_5$ is as follows.
\[
\begin{array}{c|rrrrrrr}
  \rm class&\rm1&\rm2A&\rm2B&\rm3A&\rm4A&\rm5A&\rm6A\cr
  \rm size&1&10&15&20&30&24&20 \cr
\hline
  \eta_{1}&1&1&1&1&1&1&1 \cr
  \eta_{2}&1&-1&1&1&-1&1&-1 \cr
  \eta_{3}&4&-2&0&1&0&-1&1 \cr
  \eta_{4}&4&2&0&1&0&-1&-1 \cr
  \eta_{5}&5&1&1&-1&-1&0&1\cr
  \eta_{6}&5&-1&1&-1&1&0&-1\cr
  \eta_{7}&6&0&-2&0&0&1&0\cr
\end{array}
\qquad
\begin{array}{l}
\rm 2A : (12) \\
\rm 2B : (12)(34)\\
\rm 3A : (123) \\
\rm 4A : (1234)\\
\rm 5A : (12345) \\
\rm 6A : (12)(345)
\end{array}
\]
Since we have $\widehat{\Sy_5} = \C_2$, we consider only invariants with the trivial character. 
Faithful characters satisfying the conditions on absolute values are
\begin{align*}
\phi_1 =  \eta_1 + \eta_2 + \eta_4, \quad \phi_1 \otimes \eta_2 = \eta_1 + \eta_2 + \eta_3, \\
\phi_2 = \eta_1 + \eta_5, \quad \phi_2 \otimes \eta_2 = \eta_2 + \eta_6. 
\end{align*}
Note that we have
\[
\phi_1|\A_5 = (\phi_1 \otimes \eta_2)|\A_5 = \theta_1, \quad
\phi_2|\A_5 = (\phi_2 \otimes \eta_2)|\A_5 = \theta_2
\]
for $\theta_1$ and $\theta_2$ in (\ref{theta-A5}). By extending actions of $\A_5$ in 
section \ref{SS-A5}, we can obtain invariant cubics as subspaces of $V_3(\rho_1)$ in 
(\ref{eqA5-1}) and $V_3(\rho_2)$ in (\ref{eqA5-2}). The character $\phi_1$ is afforded by 
$\pi_1 : \Sy_5 \rightarrow \GL_6(\CC)$ such that 
\begin{align*}
\pi(\sigma)[x_1 \cdots : x_5 : y] 
= [x_{\sigma(1)} : \cdots : x_{\sigma(5)} : (\mathrm{sgn} \, \sigma) y], 
\end{align*} 
and this satisfies the symplectic conditions. We have $\pi_1|\A_5 = \rho_1$ and 
$V_3(\rho_1)$ is decomposed into $(\pm1)$-eigenspaces for $(12) \in \Sy_5$ 
\begin{align} \label{eqS5-1}
V_3(\rho_1)_+ &= V_3(\pi_1) 
= \mathrm{Span} \{ y^2 s_1(x), \quad s_1(x)^3, \quad s_1(x) s_2(x), \quad s_3(x) \} 
\end{align}
and
\begin{align}
V_3(\rho_1)_- &= V_3(\pi_1 \otimes \eta_2) 
= \mathrm{Span} \{ y^3, \quad y s_1(x)^2, \quad y s_2(x) \}.
\end{align}
A general member of $V_3(\pi_1)$ defines a smooth cubic 4-fold since $y^2s_1(x) + s_3(x) = 0$ 
is smooth, but $F \in V_3(\pi_1 \otimes \eta_2)$ is reducible. 
We have $\dim \C_{\pi_1}(\Sy_5) = 1^2 + 1^2 + 1^2 = 3$ and
\[
\dim V_3(\pi_1) - \dim \C_{\pi_1}(\A_5) = 1 = 20 - r(\Sy_5).
\]
Next we consider a representation $\pi_2$ affording $\phi_2$. 
Let us consider the following involution with eigenvalues $(1,1,1,1,-1,-1)$: 
\begin{align}
a = \frac{1}{3}[a_{ij}]_{6 \times 6} \qquad 
a_{ij} = \begin{cases} -2 \ (i + j \equiv 3 \mod 6) \\ 
1 \ (i + j \not\equiv 3 \mod 6)\end{cases}. 
\end{align}
It normalizes $\rho_2(\A_5)$ in section \ref{SS-A5}, and generates $\Sy_5$ with 
$\rho_2(\A_5)$.  This involution decomposes $V_3(\rho_2)$ into   
\begin{align}\label{eqS5-2}
V_3(\rho_2)_+ &= V_3(\pi_2) 
= \mathrm{Span} \{ s_1(x)^3, \quad s_1(x)s_2(x), \quad I_1-I_2 \} 
\end{align}
and
\begin{align}
V_3(\rho_2)_- &= V_3(\pi_2 \otimes \eta_2) 
= \mathrm{Span} \{ s_1(x)^3 - 9s_1(x)s_2(x) + 18s_3(x) \}.
\end{align}
A general member of $V_3(\pi_2)$ defines a smooth cubic 4-fold since $s_1(x)s_2(x) + (I_1-I_2) = 0$ 
is smooth. On the other hand, 
\[
s_1(x)^3 - 9s_1(x)s_2(x) + 18s_3(x) = 0
\]
is singular at $[1: \cdots : 1]$. 
We have $\dim \C_{\pi_2}(\Sy_5) = 1^2 + 1^2 = 2$ and
\[
\dim V_3(\pi_2) - \dim \C_{\pi_2}(\A_5) = 1 = 20 - r(\Sy_5).
\]
\subsection{The Mathieu group $\M_9$} \label{SS-M9}
We have $\M_9 \cong 3^2.\Q_8$ and this group contains $3^2.4$ in section \ref{SS-3^2.4} 
as a subgroup of index 2. 
The Schur multiplier of $\M_9$ is $\C_3$, and we must consider a Schur cover $3.\M_9$. 
The character table of $\M_9$ itself is as follows: 
\[
\begin{array}{c|rrrrrr}
  \rm class&\rm1&\rm2A&\rm3A&\rm4A&\rm4B&\rm4C\cr
  \rm size&1&9&8&18&18&18\cr
\hline
  \chi_{1}&1&1&1&1&1&1\cr
  \chi_{2}&1&1&1&-1&1&-1\cr
  \chi_{3}&1&1&1&1&-1&-1\cr
  \chi_{4}&1&1&1&-1&-1&1\cr
  \chi_{5}&2&-2&2&0&0&0\cr
  \chi_{6}&8&0&-1&0&0&0\cr
\end{array}
\]
For any 6-dimensional representation $\rho$, we have $\Tr_{\rho}(3\A) = 6$ and 
this is not the cases.  
Next, let us consider a Schur cover $3.\M_9$ with GAP-ID (216, 88). 
Since we have $\widehat{\M_9} = \C_2 \times \C_2$, we consider only invariants 
with the trivial character.
Irreducible representaions of $3.\M_9$ with a non-trivial action of the Schur multiplier are 
given as follows.
\begin{align*}
\begin{array}{c|rrrrrrrrrrr}
  \rm class&\rm1&\rm2A&\rm3A&\rm3B&\rm3C&\rm4A&\rm4B&\rm4C&\rm6A&\rm6B&\cdots\cr
  \rm size&1&9&1&1&24&18&18&18&9&9&\cdots \cr
\hline
\chi_7    & 3 & -1 & \al & \ba & 0 &-1 &-1 & 1 & \bg & \ga \\
\chi_8    & 3 & -1 & \al & \ba & 0 & -1 & 1 &-1 & \bg & \ga \\
\chi_9    & 3 & -1 & \al & \ba & 0 &  1 &-1 &-1 & \bg & \ga  \\
\chi_{10} & 3 & -1 & \al & \ba & 0 &  1 & 1 & 1 & \bg & \ga  \\
\chi_{11} & 3 & -1 & \ba & \al & 0 & -1 &-1 & 1 & \ga & \bg   \\
\chi_{12} & 3 & -1 & \ba & \al & 0 & -1 & 1 &-1 & \ga & \bg  \\
\chi_{13} & 3 & -1 & \ba & \al & 0 &  1 &-1 &-1 & \ga & \bg \\
\chi_{14} & 3 & -1 & \ba & \al & 0 & 1 & 1 & 1 & \ga & \bg    \\
\chi_{15} & 6 &  2 & \be & \bb & 0 &  0 & 0 & 0 & \bd & \de  \\
\chi_{16} & 6 &  2 & \bb & \be & 0 & 0 &  0 & 0 & \de & \bd 
\end{array}
\\
\alpha = 3\w, \quad \beta = 6\w, \quad \gamma = -\w, \quad \delta = 2\w
\end{align*}
Classes 1, 3A and 3B form the Schur multiplier $\C_3$ which act as scalar matrices for
\[
\chi_i + \chi_j \quad (7 \leq i,j \leq 10), \qquad
\chi_k + \chi_l \quad (11 \leq k,l \leq 14). \quad 
\chi_{15}, \quad \chi_{16}. 
\]
Among these characters, conditions 
\[
\chi({\rm 4A}) = \chi({\rm 4B}) = \chi({\rm 4C}) = 0
\]
are satisfied only by $\chi_{15}$ and $\chi_{16} = \overline{\chi_{15}}$.  Representaions 
affording $\chi_{15}$ and $\overline{\chi_{15}}$ are given by extending the action of 
$H_2 \cong 3.(3^2.4)$ in section \ref{SS-3^2.4}. Let us consider a matrix 
$g = \begin{bmatrix} {\bf 0} & h_2 \\ \w h_2 & {\bf 0} \end{bmatrix}$ 
of order 4 with eigenvalues $(1,1,-1,-1,i,-i)$ where
\begin{align}
h_2 = \frac{1}{\sqrt{3}} \begin{bmatrix} \w & \w^2 & \w^2 \\ 1 & 1 & \w^2 \\ 
1 & \w^2 & 1 \end{bmatrix}.
\end{align}
The group $H_2$ together with $g$ generate a group $H_3$ with 
GAP-ID (216, 88), and this representation $\rho$ affords $\chi_{15}$ or 
$\overline{\chi_{15}}$. Moreover we have $\overline{H_3} = \delta H_3 \delta^{-1}$ for 
$\delta = \diag(1, 1, 1, i \w, i \w, i \w)$ since we have $\overline{H_2} = H_2$ 
(see section \ref{SS-3^2.4}) and
\[
\overline{g} = \delta g^{-1} \tilde{h} \delta^{-1} \in \delta H_3 \delta^{-1}, \qquad 
\delta H_2 \delta^{-1} = H_2 
\]
where $\tilde{h}$ was given in section \ref{SS-3^2.4}.  
Therefore we have $V_3(\overline{\rho}) = \delta^*V_3(\rho)$, and 
we do not have to consider $\overline{\rho}$. The space $V_3(\rho)$ consists of 
$g$-invariants in (\ref{eq3^2.4-2}):
\begin{multline} \label{eqM9}
\mathrm{Span} \{ (s_3(x) - 3u_- x_1 x_2 x_3) + (\w + i \w^2)(s_3(y) - 3u_+ y_1 y_2 y_3), \\
(i - \w^2)(I_1 - u_- I_2) + (I_3 - u_+ I_4) \}.
\end{multline}
We have $\dim \C_{\rho}(3.\M_9) = 1$ and
\[
\dim V_3(\rho) - \dim \C_{\rho}(3.\M_9) = 1 = 20 - r(\M_9).
\]
\subsection{The group $\mathbf{N_{72} \cong (S_3 \times S_3):2}$ 
with GAP-ID (72,40)}
This is a subgroup of $\Sy_6$ generated by
\[
(12), \quad (123), \quad (45), \quad (456), \quad (14)(25)(36).
\]
The Schur multiplier of $\mathrm{N}_{72}$ is $\C_2$, and we do not have to consider a
Schur cover. The character table of $\mathrm{N}_{72}$ is as follows.
\[
\begin{array}{c|rrrrrrrrr}
  \rm class&\rm1&\rm2A&\rm2B&\rm2C&\rm3A&\rm3B&\rm4A&\rm6A&\rm6B\cr
  \rm size&1&6&6&9&4&4&18&12&12\cr
\hline
  \xi_{1}&1&1&1&1&1&1&1&1&1\cr
  \xi_{2}&1&-1&-1&1&1&1&1&-1&-1\cr
  \xi_{3}&1&-1&1&1&1&1&-1&-1&1\cr
  \xi_{4}&1&1&-1&1&1&1&-1&1&-1\cr
  \xi_{5}&2&0&0&-2&2&2&0&0&0\cr
  \xi_{6}&4&-2&0&0&1&-2&0&1&0\cr
  \xi_{7}&4&0&-2&0&-2&1&0&0&1\cr
  \xi_{8}&4&0&2&0&-2&1&0&0&-1\cr
  \xi_{9}&4&2&0&0&1&-2&0&-1&0\cr
\end{array}
\qquad
\begin{array}{l}
\rm 2A = (12) \\
\rm 2B = (14)(25)(36) \\
\rm 2C = (12)(45) \\
\rm 3A = (123) \\
\rm 3B = (123)(456) \\
\rm 4A = (14)(2536) \\
\rm 6A =  (12)(456) \\
\rm 6B = (142536) \\
\end{array}
\]
Since we have $\widehat{\mathrm{N}_{72}} \cong C_2 \times \C_2$,we consider only 
invariants with the trivial character. We see that $\mathrm{N}_{72}$-invariant 
cubics are given as subspaces of $\Sy_{3,3}$-invariant cubics $V_3(\rho_1)$ 
in (\ref{eqS33}). We have
\begin{align*}
\xi_1|\Sy_{3,3} = \xi_4|\Sy_{3,3} = \eta_{11}, \quad
\xi_2|\Sy_{3,3} = \xi_3|\Sy_{3,3} = \eta_{22}, \quad
\xi_5|\Sy_{3,3} = \eta_{21} + \eta_{12}, \\
\xi_6|\Sy_{3,3} = \eta_{23} + \eta_{32}, \quad \xi_9|\Sy_{3,3} = \eta_{13} + \eta_{31}, \quad
\xi_7|\Sy_{3,3} = \xi_8|\Sy_{3,3} = \eta_{33}
\end{align*}
for $\eta_{ij}$ in section \ref{SS-S33}, and characters $\psi$ satisfying  
\begin{align*}
|\psi(2A)| = |\psi(2B)| = 2, \quad 
\psi|\Sy_{3,3} = \phi_1 = \eta_{11} + \eta_{22} + \eta_{13} + \eta_{31}
\end{align*}
are only
\[
\psi = \xi_1 + \xi_3 + \xi_9, \quad \psi \otimes \xi_2 = \xi_2 + \xi_4 + \xi_9. 
\]
A representation $\pi$ affording the character $\psi$ is given by extending 
$\rho_1$ in secton \ref{SS-S33}. Let us consider a matrix 
$c = \begin{bmatrix} C_1 & C_2 \\ C_2 & C_1 \end{bmatrix}$ of order 2 with 
eigenvalues $(1,1,1,1,-1,-1)$, where 
\[
C_1 = \begin{bmatrix} 1 & 1& 1 \\ 1 & 1 & 1 \\ 1 & 1 & 1 \end{bmatrix}, \quad
C_2 = \begin{bmatrix} -2 & 1& 1 \\ 1 & -2 & 1 \\ 1 & 1 & -2 \end{bmatrix}. 
\]
The group $\rho_1(\Sy_{3,3})$ and $c$ generate a group with GAP-ID (72, 40), 
and this gives $\pi$. By the involution $c$, the space $V_3(\rho_1)$ is 
decomposed into the $(+1)$-eigenspace
\begin{multline} \label{eqN72}
V_3(\pi) = \mathrm{Span} \{ s_1(y)^3, \ s_1(x)^2 s_1(y), \ 
s_2(y)s_1(y) + s_2(x) s_1(y), \\ 9s_2(y)s_1(y) -9s_3(y) + p(x) \} 
\end{multline}
and the $(-1)$-eigenspace
\begin{multline}
V_3(\pi, \xi_2) = V_3(\pi \otimes \xi_2) = \mathrm{Span} 
\{ 2 s_1(y)^3 -9 s_2(y)s_1(y) - 9 s_3(y) + p(x) ,\\
s_1(y)^3 -3 s_2(y)s_1(y) - s_1(x)^2 s_1(y) + 3 s_2(x) s_1(y) \}. 
\end{multline}
A general member of $V_3(\pi)$ defines a smooth cubic fourfold since 
\[
s_1(x)^2 s_1(y) + (9s_2(y)s_1(y) -9s_3(y) + p(x)) = 0
\]
is smooth. For $F \in V_3(\pi \otimes \chi_2)$, a cubic fourfold $F = 0$ is singular 
along $x_1 = x_2 = x_3,\ y_1 = y_2 = y_3$.
We have $\dim \C_{\pi}(\Sy_{3,3}:2) = 1^2 + 1^2 + 1^2 = 3$ and
\[
\dim V_3(\pi) - \dim \C_{\pi}(\Sy_{3,3}:2) = 1 = 20 - r(\Sy_{3,3}:2).
\]
\subsection{The group $\mathbf{T_{48} \cong Q_8:S_3}$ with GAP-ID (48,29)}
The Schur multiplier of $\mathrm{T}_{48}$ is 1, and hence we have 
$(\mathrm{T}_{48})^* = \mathrm{T}_{48}$. The character table of $\mathrm{T}_{48}$ is as follows.
\[
\begin{array}{c|rrrrrrrr}
  \rm class&\rm1&\rm2A&\rm2B&\rm3A&\rm4A&\rm6A&\rm8A&\rm8B\cr
  \rm size&1&1&12&8&6&8&6&6\cr
\hline
  \chi_{1}&1&1&1&1&1&1&1&1\cr
  \chi_{2}&1&1&-1&1&1&1&-1&-1\cr
  \chi_{3}&2&2&0&-1&2&-1&0&0\cr
  \chi_{4}&2&-2&0&-1&0&1&-\sqrt{-2}&\sqrt{-2}\cr
  \chi_{5}&2&-2&0&-1&0&1&\sqrt{-2}&-\sqrt{-2}\cr
  \chi_{6}&3&3&-1&0&-1&0&1&1\cr
  \chi_{7}&3&3&1&0&-1&0&-1&-1\cr
  \chi_{8}&4&-4&0&1&0&-1&0&0\cr
\end{array}
\]
Since we have $\widehat{\mathrm{T}_{48}} = \C_2$, we consider only invariants 
with the trivial character. Characters $\theta$ of degree 6 satisfying
\[
|\theta(2A)| = |\theta(2B)| = 2, \quad  |\theta(4A)| = 0, \quad 
|\theta(8A)| = |\theta(8B)| = \sqrt{2}
\]
are only
\begin{align*}
\theta_1 = \chi_1 + \chi_4 + \chi_7, \quad 
\theta_1 \otimes \chi_2 = \chi_2 + \chi_5 + \chi_6, \\
\overline{\theta_1} = \chi_1 + \chi_5 + \chi_7, \quad 
\overline{\theta_1 \otimes \chi_2} = \chi_2 + \chi_4 + \chi_7. 
\end{align*}
Since we have Molien series
\begin{align*}
\mu(\theta_1, t) = 1 + t + 2t^2 + 4t^3 + \cdots, \qquad
\mu(\theta_1 \otimes \chi_2, t) = 1 + 2t^2 + 7t^4 + \cdots,
\end{align*}
we consider only $\theta_1$ and $\overline{\theta_1}$. The group $T_{48}$ 
contains $\Q_8$ as a normal subgroup, and a representation $\pi$ affording 
$\theta_1$ is given by extending the representation $\rho$ in section \ref{SS-Q8}. 
Let us consider two matrices  
\begin{align*}
n_1 = 
\begin{bmatrix} 
1 & 0 & 0 & 0 \\ 0 & 0 & 1 & 0 \\ 0 & 1 & 0 & 0 \\ 0 & 0 & 0 & 1
\end{bmatrix}
\oplus \frac{1}{\sqrt{2}}
\begin{bmatrix} 1 & -i \\ i & -1 \end{bmatrix}, \quad
n_2 = \begin{bmatrix} 1 & 0 \\ 0 & 1 \end{bmatrix}
\oplus \begin{bmatrix} 0 & 1 \\ 1 & 0 \end{bmatrix}
\oplus \begin{bmatrix} 0 & \z_8^3 \\ \z_8^5 & 0 \end{bmatrix}
\end{align*}
of order 2 with eigenvalues $(1,1,1,1,-1,-1)$. They normalize $\rho(\Q_8)$, 
and $n_1, n_2$ together with $\rho(\Q_8)$ generate a group 
$H \cong \Q_8:\Sy_3$ with GAP-ID (48,29) (Conjugations by $n_1$ and $n_2$ permute 
conjugacy classes 
\[
\{ \rho(I),\ \rho(-I) \}, \quad \{ \rho(J),\ \rho(-J) \}, \quad 
\{ \rho(K),\ \rho(-K) \}
\]
of $\rho(\Q_8)$, and they generate $\mathrm{Out}(\Q_8) \cong \Sy_3$). 
Moreover we have $\overline{H} = \delta H \delta^{-1}$ for 
$\delta = \diag(1,1,-1,-1) \oplus \begin{bmatrix} 0 & 1 \\ 1 & 0 \end{bmatrix}$ 
since
\begin{align*}
\overline{\rho(I)} = \delta \rho(I)^3 \delta^{-1}, \quad
 \overline{\rho(J)} = \delta \rho(J) \delta^{-1}, \quad
\overline{n_1} = \delta (n_1 \rho(IJ)) \delta^{-1}, \quad
\overline{n_2} = \delta (n_2) \delta^{-1}.
\end{align*}
Therefore we have $V_3(\overline{\pi}) = \delta^*V_3(\pi)$ and we do not have to 
consider $V_3(\overline{\pi})$.  The space $V_3(\pi)$ is a subspace of (\ref{eqQ8}):
\begin{multline} \label{eqT48}
\mathrm{Span} \{ x_1^3, \quad x_1(x_2^2 + x_3^2 + x_4^2), \quad x_2 x_3 x_4, \\
i x_3 (x_5^2 + x_6^2)  + 2 x_2 x_5 x_6 + x_4 (x_5^2 - x_6^2) \}.
\end{multline}
A general member of $V_3(\pi)$ is smooth since sum of the above four generator define 
a smooth cubic fourfold.
We have $\dim \C_{\pi}(\mathrm{T}_{48}) = 1^2 + 1^2 + 1^2 = 3$ and
\[
\dim V_3(\pi) - \dim \C_{\pi}(\mathrm{T}_{48}) = 1 = 20 - r(\mathrm{T}_{48}).
\]

%% file: maximal.tex
\section{Maximal groups}
We have six maximal groups:
\begin{align*}
3^4:\A_6, \quad \A_7, \quad 3^{1+4}_+:2.2^2, \quad \M_{10}, \quad 
\mathrm{L}_2(11), \quad \A_{3,5}.
\end{align*}
There exists the unique cubic fourfold $X$ such that $\Aut^s(X)= G$ except for 
$G = \A_7$ and $\M_{10}$. However, there exist two non-isomorphic cubic fourfolds $X$ 
such that $\Aut^s(X)$ is $\A_7$ or $\M_{10}$.
Corresponding eight cubic fourfolds were identified in \cite{HM19}, \cite{LZ22} and 
\cite{YYZ23}.
\\ \\ 
(i) \ For $G = 3^4:\A_6$, the corresponding cubic fourfold is the Fermat cubic
\begin{align}
x_1^3 + x_2^3 + x_3^3 + x_4^3 + x_5^3 + x_6^3 = 0.
\end{align}
(ii) \ For $G = \A_7$, corresponding two cubic fourfolds are (\ref{Clebsch-Segre}) 
and (\ref{3.A7-invariant}). 
\\ \\
(iii) \ For $G = 3^{1+4}_+:2.2^2$, the corresponding cubic fourfold is given by 
\[
F(x_1,x_2,x_3) = F(x_4, x_5, x_6)
\]
 where  
\begin{align}
F(x,y,z) = x^3 + y^3 + z^3 - 3(1 + \sqrt{3})xyz.
\end{align}
Note that $F(x,y,z) = 0$ is an elliptic curve whose $j$-invariant is 
1728 (see \cite{Ko22}), and we may replace $F$ by 
\begin{align}
F(x,y,z) = y^2z - x^3 - xz^2.  
\end{align}
In particular, this cubic fourfold is defined over $\QQ$. 
\\ \\
(iv) \ For $G = \mathrm{L}_2(11) = \PSL_2(\mathbb{F}_{11})$, the corresponding cubic fourfold is given by
\begin{align}
x_6^3 = x_1^2x_2 + x_2^2 x_3 + x_3^2 x_4 + x_4^2 x_5 + x_5^2 x_1. 
\end{align}
(v) \ For $G = \A_{3,5} = (\Sy_3 \times \Sy_8) \cap \A_8$, the corresponding cubic fourfold is given by 
\begin{align}
\begin{cases} x_1^3 + x_2^3 + x_3^3 + x_4^3 + x_5^3 + x_6^3 + x_7^3 + x_8^3 = 0 \\
x_1 + x_2 + x_ 3 = 0 \\
x_4 + x_5 + x_ 6 + x_7 + x_8 = 0 
\end{cases}
\end{align}
(vi) \ For $G = \M_{10}$, we give corresponding cubic fourfolds $X_+$ and $X_-$ 
defined over $\QQ(\sqrt{6})$. 
The Mathieu group $\M_{10}$ contains $\A_6$ as a subgroup of index 2, and $X_{\pm}$ 
belongs to (\ref{eqA6-1}) or (\ref{eqA6-2}).  As in the case of $\M_9 \subset \M_{10}$, 
there is not a faithful 6-dimensional representation of $\M_{10}$ (see section \ref{SS-M9}). 
Therefore both of $X_+$ and $X_-$ belong to the family (\ref{eqA6-2}). 
\\ \indent
Let us consider a matrix 
\begin{align*}
r = \frac{1}{3 \sqrt{6}}
\begin{bmatrix}
2\w+3 & -\w & -2(\w+1) & -3\w-1 & \w-2 & -7 \\
-\w-3 & -\w & 7(\w+1) & 3\w+2 & 1-2\w & 2 \\
2\w & -4\w & \w-2 & 3\w+2 & \w-2 & 3\w+2\\
-7\w & -4\w & 1-2\w & -3\w-1 & 1-2\w & -3\w-1 \\
2\w+3 & -\w & \w-2 & -7 & -2(\w+1) & -3\w-1 \\
-\w-3 & -\w & 1-2\w & 2 & 7(\w+1) & 3\w+2 
\end{bmatrix}
\end{align*}
of order 4 with eigenvalues $(1,1,-1,-1,i,-i)$. The group $H \cong 3.\A_6$ in section \ref{SS-A6} 
is normalized by $r$, and $H$ together with $r$ generate $3.\M_{10}$. By the action of $r$, 
the space $V_3(\rho) = \CC f(x) \oplus \CC g(x)$ is decomposed into $(\pm)$-eigenspaces
$V_3(\rho)_{\pm} = \CC F_{\pm}(x)$ where
\begin{align}
F_+(x) = f(x) - (7 + 3 \sqrt{6}) g(x), \qquad F_-(x) = f(x) - (7 - 3 \sqrt{6}) g(x).
\end{align}
Cubic fourfolds $X_{\pm}$ defined by $F_{\pm}(x) = 0$ are smooth, defined over $\QQ(\sqrt{6})$ 
and $\Aut^s(X_{\pm}) = \M_{10}$. Since we have 
\[
F_+(rx) = F_+(x), \qquad F(rx) = - F_-(x)
\]
these are not projectively equivalent, that is, not isomorphic. Moreover $X_+$ is 
 $\mathrm{Gal}(\QQ(\sqrt{6})/\QQ)$-conjugate to $X_-$.